\documentclass[12pt]{amsart}
\usepackage{amsfonts}
\usepackage{amssymb}
\usepackage{latexsym}
\usepackage{epsfig}

\newtheorem{theorem}{Theorem}[section]
\newtheorem{corollary}[theorem]{Corollary}
\newtheorem{lemma}[theorem]{Lemma}
\newtheorem{proposition}[theorem]{Proposition}

\newtheorem{problem}[theorem]{Problem}

\theoremstyle{definition}
\newtheorem{definition}[theorem]{Definition}
\newtheorem{remark}[theorem]{Remark}
\newtheorem{example}[theorem]{Example}

\begin{document}

\title{An Algebraic Perspective of Group Relaxations}
\author{Rekha R. Thomas}
\address{Department of Mathematics, University of Washington,
Box 354350, Seattle, WA 98195}
\email{thomas@math.washington.edu}

\date{\today}

\maketitle

\tableofcontents 

\section{Introduction}
{\em Group relaxations} of integer programs were introduced by Ralph
Gomory in the 1960s \cite{Gom65}, \cite{Gom69}. Given a general 
integer program of the form 
\begin{equation} \label{generalip}
minimize \,\,\{ c \cdot x
\,\, : \,\, Ax = b, \,\, x \geq 0, \,\,integer \},
\end{equation}
its group relaxation is obtained by dropping non-negativity
restrictions on all the basic variables in the optimal solution of its
linear relaxation. In this paper, we survey recent results on group
relaxations obtained from the algebraic study of integer programming
using {\em Gr\"obner bases} of {\em toric ideals} \cite{Stu}. No
knowledge of these methods is assumed, and the exposition is
self-contained and hopefully accessible to a person familiar with the 
traditional methods of integer programming. Periodic comments on the
algebraic origins, motivations and counterparts of many of the
described results --- which the reader may pursue if desired --- offer
a more complete picture of the theory.

For the sake of brevity, we will bypass a detailed account of the
classical theory of group relaxations. A short expository account can
be found in \cite[\S 24.2]{Sch}, and a detailed set of lecture notes 
on this topic in \cite{Joh}. We give a brief synopsis of the
essentials based on the recent survey article by Aardal et. al
\cite{AWW} and refer the reader to any of the above sources for
further details and references on the classical theory of group
relaxations.

Assuming that all data in (\ref{generalip}) are integral and that
$A_B$ is the optimal basis of the linear relaxation of
(\ref{generalip}), Gomory's group relaxation of (\ref{generalip}) is
the problem
\begin{equation} \label{generalgr}
minimize \,\,\{\tilde c \cdot x_N \,\,
  : \,\, A_B^{-1}A_Nx_N \equiv A_B^{-1}b \,\, (mod \,\,1),\,\,x_N \geq 0,
  \,\,integer \}.
\end{equation}
Here $B$ and $N$ are the index sets for the basic and non-basic
columns of $A$ corresponding to the optimal solution of the linear
relaxation of (\ref{generalip}). The vector $x_N$ denotes the
non-basic variables and the cost vector $\tilde c = c_N -
c_BA_B^{-1}A_N$ where $c= (c_B,c_N)$ is partitioned according to $B$
and $N$. The notation $A_B^{-1}A_Nx_N \equiv A_B^{-1}b \,\, (mod
\,\,1)$ indicates that $A_B^{-1}A_Nx_N- A_B^{-1}b$ is a vector of
integers. Problem (\ref{generalgr}) is called a ``group relaxation''
of (\ref{generalip}) since it can be written in the canonical form
\begin{equation} \label{generalgp}
minimize \,\,\{\tilde c \cdot x_N \,\,: \,\, \sum_{j \in N} g_jx_j
\equiv g_0 \,\,(mod \,\, G),\,\, x_N \geq 0,\,\,integer\}
\end{equation}
where $G$ is a finite abelian group and $g_j \in G$. Problem
(\ref{generalgp}) can be viewed as a shortest path problem in a graph
on $|G|$ nodes which immediately furnishes algorithms for solving it.
Once the optimal solution $x_N^*$ of (\ref{generalgr}) is found, it
can be uniquely {\em lifted} to a vector $x^* = (x_B^*,x_N^*) \in
\mathbb Z^n$ such that $Ax^* = b$. If $x_B^* \geq 0$ then $x^*$ is the
optimal solution of (\ref{generalip}).  Otherwise, $c \cdot x^*$ is a
lower bound for the optimal value of (\ref{generalip}).  Several
strategies are possible when the group relaxation fails to solve the
integer program.  See \cite{BeSh},\cite{GNS}, \cite{NW} and
\cite{Wol73} for work in this direction. A particular idea due to
Wolsey \cite{Wol71} that is very relevant for this paper is to
consider the {\em extended} group relaxations of (\ref{generalip}).
These are all the possible group relaxations of (\ref{generalip})
obtained by dropping non-negativity restrictions on all possible
subsets of the basic variables $x_B$ in the optimum of the linear
relaxation of (\ref{generalip}). Gomory's group relaxation
(\ref{generalgr}) of (\ref{generalip}) and (\ref{generalip}) itself
are therefore among these extended group relaxations. If
(\ref{generalgr}) does not solve (\ref{generalip}), then one could
resort to other extended relaxations to solve the problem. At least
one of these extended group relaxations (in the worst case
(\ref{generalip}) itself) is guaranteed to solve the integer program 
(\ref{generalip}).

The convex hull of the feasible solutions to (\ref{generalgr}) is
called the {\em corner polyhedron} \cite{Gom67}.  A major focus of
Gomory and others who worked on group relaxations was to understand
the polyhedral structure of the corner polyhedron. This was achieved
via the {\em master polyhedron} of the group $G$ \cite{Gom69} which is
the convex hull of the set of points
$$
\{z \,\, : \,\,\sum_{g \in G} gz_g
\equiv g_0 \,\,(mod \,\, G),\,\, z \geq 0,\,\,integer\}. 
$$
Facet-defining inequalities for the master polyhedron provide facet
inequalities of the corner polyhedron \cite{Gom69}. As remarked in
\cite{AWW}, this landmark paper \cite{Gom69} introduced several of the
now standard ideas in polyhedral combinatorics like projection onto
faces, subadditivity, master polytopes, using automorphisms to
generate one facet from another, lifting techniques and so on. See
\cite{GJ} for further results on generating facet inequalities.

In the algebraic approach to integer programming, one considers the
entire family of integer programs of the form (\ref{generalip}) as the
right hand side vector $b$ varies. Definition~\ref{grouprels} defines
a set of group relaxations for each program in this family.  Each
relaxation is indexed by a face of a simplicial complex called a {\em
  regular triangulation} (Definition~\ref{regtriang}). This complex
encodes all the optimal bases of the linear programs arising from the
coefficient matrix $A$ and cost vector $c$ (Lemma~\ref{optbases}). The
main result of Section~2 is Theorem~\ref{bounded} which states that
the group relaxations in Definition~\ref{grouprels} are precisely all
the bounded group relaxations of all programs in the family. In
particular, they include all the extended group relaxations of all
programs in the family and typically contain more relaxations for each
program. This theorem is proved via a particular reformulation of
group relaxations which is crucial for the rest of the paper. This and
other reformulations are described in Section~2.

The most useful group relaxations of an integer program are the
``least strict'' ones among all those that solve the program. By this
we mean that any further relaxation of non-negativity restrictions
will result in group relaxations that do not solve the problem. The
faces of the regular triangulation indexing all these special
relaxations for all programs in the family are called the {\em
  associated sets} of the family (Definition~\ref{associatedset}).  In
Section~3 we develop tools to study associated sets. This leads to
Theorem~\ref{stdpairs-equivs} which characterizes associated sets in
terms of {\em standard pairs} and {\em standard
  polytopes}. Theorem~\ref{cover=solve} shows that one can ``read
off'' the ``least strict'' group relaxations that solve a given
integer program in the family from these standard pairs.

The results in Section~3 lead to an important invariant of the family
of integer programs being studied called its {\em arithmetic degree}.
In Section~4 we discuss the relevance of this invariant and give a
bound for it based on a result of Ravi Kannan (Theorem~\ref{kannan}).
His result builds a bridge between our methods and those of Kannan,
Lenstra, Lovasz, Scarf and others that use geometry of numbers in
integer programming.

Section~5 examines the structure of the poset of associated sets.  The
main result in this section is the {\em chain theorem}
(Theorem~\ref{chain-theorem}) which shows that associated sets occur
in saturated chains. Theorem~\ref{length} bounds the length of a
maximal chain.

In Section~6 we define a particular family of integer programs called
a {\em Gomory family}, for which all associated sets are maximal faces
of the regular triangulation. Theorem~\ref{gomoryfam-theorem} gives
several characterizations of Gomory families. We show that this notion
generalizes the classical notion of {\em total dual integrality} in
integer programming \cite[\S 22]{Sch}. We conclude in Section~7 with
constructions of Gomory families from matrices whose columns form a
Hilbert basis. In particular, we recast the existence of a Gomory
family as a {\em Hilbert cover} problem. This builds a connection to
the work of Seb{\"o} \cite{Seb}, Bruns \& Gubeladze \cite{BG} and
Firla \& Ziegler \cite{FZ} on {\em Hilbert partitions} and {\em
  covers} of polyhedral cones. We describe the notions of {\em super}
and $\Delta$-{\em normality} both of which give rise to Gomory
families (Theorems~\ref{specialinitial} and \ref{smalld}).

The majority of the material in this paper is a translation of
algebraic results from \cite{HT01}, \cite{HT2}, \cite{HT1}, \cite[\S 8
and \S 12.D]{Stu}, \cite{STV} and \cite{SWZ}. The translation has
sometimes required new definitions and proofs. Kannan's theorem in
Section~4 has not appeared elsewhere. 

We will use the letter $\mathbb N$ to denote the set of non-negative
integers, $\mathbb R$ to denote the real numbers and $\mathbb Z$ for
the integers. The symbol $P \subseteq Q$ denotes that $P$ is a subset
of $Q$, possibly equal to $Q$, while $P \subset Q$ denotes that $P$ is
a proper subset of $Q$.

\section{Group Relaxations}
Throughout this paper, we fix a matrix $A \in \mathbb Z^{d \times n}$
of rank $d$, a cost vector $c \in \mathbb Z^n$ and consider the family
$IP_{A,c}$ of all integer programs
$$IP_{A,c}(b) := minimize \,\,\{ c \cdot x \,\, : \,\, Ax = b, \,\, x
\in \mathbb N^n \}$$
as $b$ varies in the semigroup $\mathbb N A := \{
Au \, : \, u \in \mathbb N^n \} \subseteq \mathbb Z^d$. This family is
precisely the set of all feasible integer programs 
with coefficient matrix $A$ and cost vector $c$. The semigroup
$\mathbb N A$ lies in the intersection of the $d$-dimensional
polyhedral cone $cone(A) := \{ Au \,\, : \,\, u \geq 0 \} \subseteq
\mathbb R^d$ and the $d$-dimensional lattice $\mathbb Z A := \{ Au
\,\, : \,\, u \in \mathbb Z^n \} \subseteq \mathbb Z^d$. For
simplicity, we will assume that $cone(A)$ is pointed and that $\{u \in
\mathbb R^n \,\,:\,\,Au = 0\}$, the kernel of $A$, intersects the
non-negative orthant of $\mathbb R^n$ only at the origin. This
guarantees that all programs in $IP_{A,c}$ are bounded.  In addition,
the cost vector $c$ will be assumed to be {\em generic} in the sense
that each program in $IP_{A,c}$ has a unique optimal solution.

The {\em linear relaxation} of $IP_{A,c}(b)$ is the linear program
$$LP_{A,c}(b) := minimize \,\,\{\,\, c \cdot x \,\, : \,\, Ax = b,
\,\, x \geq 0 \}.$$
We denote by $LP_{A,c}$ the family of all linear
programs of the form $LP_{A,c}(b)$ as $b$ varies in $cone(A)$. These
are all the feasible linear programs with coefficient matrix $A$ and
cost vector $c$. Since all data are integral and all programs in
$IP_{A,c}$ are bounded, all programs in $LP_{A,c}$ are bounded as
well.

In the classical definitions of group relaxations of $IP_{A,c}(b)$,
one assumes knowledge of the optimal basis of the linear relaxation
$LP_{A,c}(b)$. In the algebraic set up, we define group
relaxations for all members of $IP_{A,c}$ at one shot and,
analogously to the classical setting, assume that the optimal 
bases of all programs in $LP_{A,c}$ are known. This information is
carried by a {\em polyhedral complex} called the {\em regular
triangulation} of $cone(A)$ with respect to $c$.

A polyhedral complex $\Delta$ is a collection of
polyhedra called {\em cells (or faces)} of $\Delta$ such that: \\
(i) every face of a cell of $\Delta$ is again a cell of $\Delta$ and, \\
(ii) the intersection of any two cells of $\Delta$ is a common face of
both. \\ 
The set-theoretic union of the cells of $\Delta$ is called the {\em
support} of $\Delta$. If $\Delta$ is not empty, then the empty set is a
cell of $\Delta$ since it is a face of every polyhedron. If all the
faces of $\Delta$ are cones, we call $\Delta$ a {\em cone
  complex}.

For $\sigma \subseteq \{1, \ldots, n\}$, let $A_{\sigma}$ be the
submatrix of $A$ whose columns are indexed by $\sigma$, and let
$cone(A_{\sigma})$ denote the cone generated by the columns of
$A_{\sigma}$. The {\em regular subdivision} $\Delta_c$ of $cone(A)$ is
a cone complex with support $cone(A)$ defined as follows.

\begin{definition} \label{regtriang}
  For $\sigma \subseteq \{1, \ldots, n \}$, $cone(A_{\sigma})$ is a
  face of the {\em regular subdivision} $\Delta_c$ of $cone(A)$ if and
  only if there exists a vector $y \in \mathbb R^d$ such that $y \cdot
  a_j = c_j$ for all $j \in \sigma$ and $y \cdot a_j < c_j$ for all $j
  \not \in \sigma$.
\end{definition}

The regular subdivision $\Delta_c$ can be constructed geometrically as
follows. Consider the cone in $\mathbb R^{d+1}$ generated by the
lifted vectors $(a_i^t,c_i) \in \mathbb R^{d+1}$ where $a_i$ is the
$i$th column of $A$ and $c_i$ is the $i$th component of $c$. The
{\em lower} facets of this lifted cone are all those facets whose
normal vectors have a negative $(d+1)$th component. Projecting these
lower facets back onto $cone(A)$ induces the regular subdivision
$\Delta_c$ of $cone(A)$ (see \cite{BFS}). Note that if the columns of
$A$ span an affine hyperplane in $\mathbb R^d$, then $\Delta_c$ can
also be seen as a subdivision of $conv(A)$, the $(d-1)$-dimensional
convex hull of the columns of $A$.

The genericity assumption on $c$ implies that $\Delta_c$ is in fact a
triangulation of $cone(A)$ (see \cite{ST}). We call $\Delta_c$ the
{\em regular triangulation} of $cone(A)$ with respect to $c$.  For
brevity, we may also refer to $\Delta_c$ as the regular triangulation
of $A$ with respect to $c$. Using $\sigma$ to label
$cone(A_{\sigma})$, $\Delta_c$ is usually denoted as a set of subsets
of $\{1,\ldots, n \}$. Since $\Delta_c$ is a complex of simplicial
cones, it suffices to list just the maximal elements (with respect to
inclusion) in this set of sets. By definition, every one dimensional
face of $\Delta_c$ is of the form $cone(a_i)$ for some column $a_i$ of
$A$. However, not all cones of the form $cone(a_i)$, $a_i$ a column of
$A$, need appear as a one dimensional cell of $\Delta_c$.

\begin{example} \label{exs_regtriangs}
(i) Let $A = \left[ \begin{array}{cccc} 1 & 1 & 1 & 1 \\ 0 & 1 & 2 & 3
\end{array} \right ]$ and $c = (1,0,0,1)$. The four columns of $A$ 
are the four dark points in Figure~\ref{f1} labeled by their column
indices $1, \ldots, 4$. Figure~\ref{f1} (a) shows the cone generated
by the lifted vectors $(a_i^t,c_i) \in \mathbb R^3$. The rays
generated by the lifted vectors have the same labels as the points
that were lifted. Projecting the lower facets of this lifted cone back
onto $cone(A)$, we get the regular triangulation $\Delta_c$ of
$cone(A)$ shown in Figure~\ref{f1} (b). The same triangulation is
shown as a triangulation of $conv(A)$ in Figure~\ref{f1} (c).
The faces of the triangulation $\Delta_c$ are $\{1,2\},\{2,3\},
\{3,4\},\{1\},\{2\},\{3\},\{4\}$ and $\emptyset$. Using only the
maximal faces, we may write $\Delta_c = \{\{1,2\}, \{2,3\},
\{3,4\}\}$.

\begin{figure}
\epsfig{file=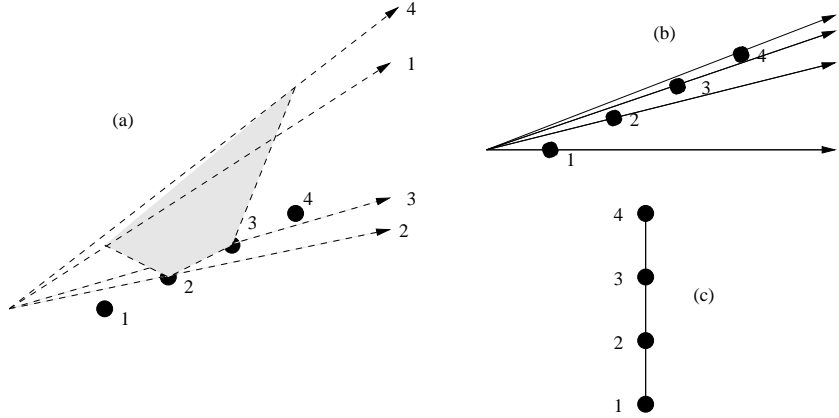, height=5.5cm}
\caption{Regular triangulation $\Delta_c$ for $c = (1,0,0,1)$
  (Example~\ref{exs_regtriangs} (i)).}
\label{f1}
\end{figure}

\noindent (ii) For the $A$ in (i), $cone(A)$ has four distinct regular
triangulations as $c$ varies. For instance, the cost vector $c' =
(0,1,0,1)$ induces the regular triangulation $\Delta_{c'} = \{
\{1,3\},\{3,4\} \}$ shown in Figure~\ref{f2} (b) and (c). Notice that
$\{2\}$ is not a face of $\Delta_{c'}$.

\begin{figure}
\epsfig{file=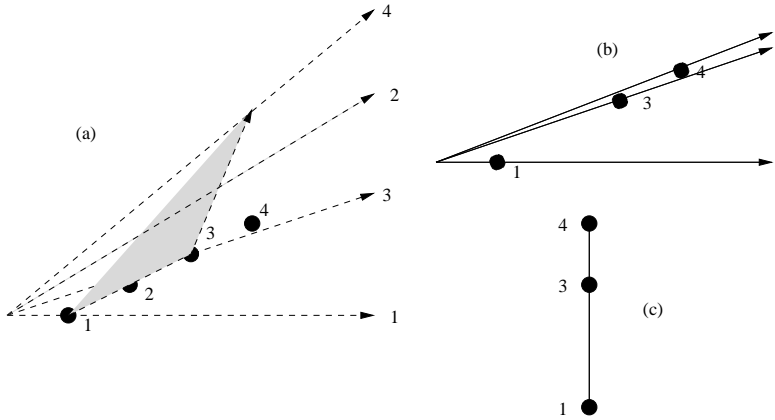, height=5.5cm}
\caption{Regular triangulation $\Delta_{c'}$ for $c' = (0,1,0,1)$
  (Example~\ref{exs_regtriangs} (ii)).}
\label{f2}
\end{figure}

\noindent (iii) If $A = \left[ \begin{array}{cccc} 1 & 3 & 2 & 1 \\ 0 &
1 & 2 & 3 \end{array} \right ]$ and $c = (1,0,0,1)$, then 
$\Delta_c = \{\{1,2\}, \{2,3\}, \{3,4\}\}.$ However, in this case,
$\Delta_c$ can only be seen as a triangulation of $cone(A)$ and not of
$conv(A)$.
\qed
\end{example}

For a vector $x \in \mathbb R^n$, let $supp(x) := 
\{i : x_i \neq 0 \}$ denote the {\em support} of $x$. The
significance of regular triangulations for linear programming is
summarized in the following proposition.

\begin{proposition} \cite[Lemma 1.4]{ST} \label{optbases}
An optimal solution of $LP_{A,c}(b)$ is any feasible solution 
$x^{\ast}$ such that $supp(x^{\ast}) = \tau$ where $\tau$ is the 
smallest face of the regular triangulation $\Delta_c$ such that $b \in
cone(A_{\tau})$. 
\end{proposition}

Proposition~\ref{optbases} implies that $\sigma \subseteq \{1, \ldots,
n\}$ is a maximal face of $\Delta_c$ if and only if $A_{\sigma}$ is an
optimal basis for all $LP_{A,c}(b)$ with $b$ in $cone(A_{\sigma})$.
For instance, in Example~\ref{exs_regtriangs} (i), 
if $b = (4,1)^t$ then the optimal basis of $LP_{A,c}(b)$ is
$[a_1,a_2]$ where as if $b = (2,2)^t$, then the optimal solution of
$LP_{A,c}(b)$ is degenerate and either $[a_1,a_2]$ or $[a_2,a_3]$
could be the optimal basis of the linear program. (Recall that $a_i$
is the $i$th column of $A$.)  All programs in $LP_{A,c}$ have one of
$[a_1,a_2],\,[a_2,a_3]$ or $[a_3,a_4]$ as its optimal basis.

Given a polyhedron $P \subset \mathbb R^n$ and a face $F$ of $P$, the
{\em normal cone} of $F$ at $P$ is the cone $N_P(F) := \{ \omega \in
\mathbb R^n \, : \, \omega \cdot x' \geq \omega \cdot x, \,\,{\text
  for \,\,all} \,
x' \in F \text{ and } x \in P\}$. The normal cones of all faces of $P$
form a cone complex in $\mathbb R^n$ called the {\it normal fan} of
$P$.

\begin{proposition} \label{normalfan}
The regular triangulation $\Delta_c$ of $cone(A)$ is the normal fan of
the polyhedron $P_c := \{ y \in \mathbb R^d \, \, : yA \leq c \}$.
\end{proposition}

\begin{proof}
  The polyhedron $P_c$ is the feasible region of $ maximize\,\,\{ y
  \cdot b \,\, : \,\, y A \leq c, \,\, y \in \mathbb R^d\}$, the dual
  program to $LP_{A,c}(b)$. The support of the normal fan of $P_c$ is
  $cone(A)$, since this is the polar cone of the recession cone $\{y
  \in \mathbb R^d : yA \leq 0 \}$ of $P_c$.  Suppose $b$ is any vector
  in the interior of a maximal face $cone(A_{\sigma})$ of $\Delta_c$.
  Then by Proposition~\ref{optbases}, $LP_{A,c}(b)$ has an optimal
  solution $x^{\ast}$ with support $\sigma$. By complementary
  slackness, the optimal solution $y$ to the dual of $LP_{A,c}(b)$
  satisfies $y \cdot a_j = c_j$ for all $j \in \sigma$ and $y \cdot
  a_j \leq c_j$ otherwise. Since $\sigma$ is a maximal face of
  $\Delta_c$, $y \cdot a_j < c_j$ for all $j \not \in \sigma$. Thus
  $y$ is unique, and $cone(A_{\sigma})$ is contained in the normal
  cone of $P_c$ at the vertex $y$. If $b$ lies in the interior of
  another maximal face $cone(A_{\tau})$ then $y'$, (the dual optimal
  solution to $LP_{A,c}(b)$) satisfies $y' \cdot A_{\tau} = c_{\tau}$
  and $y' \cdot A_{\bar \tau} < c_{\bar \tau}$ where $\tau \neq
  \sigma$. As a result, $y'$ is distinct from $y$, and each maximal
  cone in $\Delta_c$ lies in a distinct maximal cone in the normal fan
  of $P_c$. Since $\Delta_c$ and the normal fan of $P_c$ are both cone
  complexes with the same support, they must therefore coincide.
\end{proof}

\begin{figure}
\epsfig{file=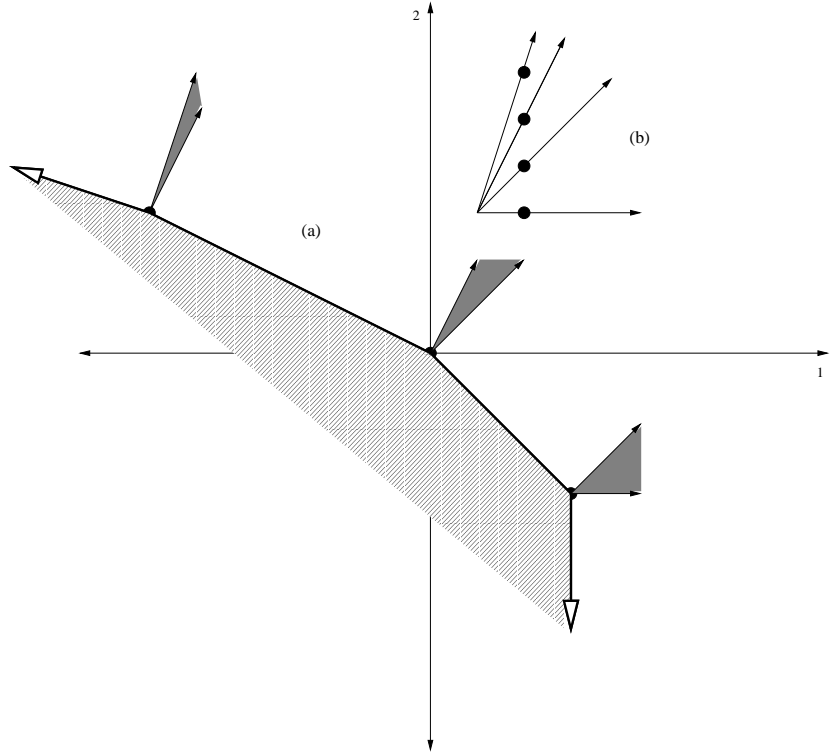, height=10cm}
\caption{The polyhedron $P_c$ and its normal fan for 
  Example~\ref{exs_regtriangs} (i).}
\label{f3}
\end{figure}

\noindent{\bf Example~\ref{exs_regtriangs} continued.}\\
Figure~\ref{f3} (a) shows the polyhedron $P_c$ for 
Example~\ref{exs_regtriangs} (i) with all its normal cones. The normal
fan of $P_c$ is drawn in Figure~\ref{f3} (b). Compare this fan with
that in Figure~\ref{f1} (b). \qed \\

\begin{corollary}
The polyhedron $P_c$ is simple if and only if the regular subdivision
$\Delta_c$ is a triangulation of $cone(A)$. 
\end{corollary}

Regular triangulations were introduced by Gel'fand, Kapranov and
Zelevinsky \cite{GKZ} and have various applications. They have played
a central role in the algebraic study of integer programming
(\cite{Stu}, \cite{ST}), and we use them now to define {\em group
  relaxations} of $IP_{A,c}(b)$.

A subset $\tau$ of $\{1, \ldots, n\}$ partitions $x = (x_1, \ldots,
x_n)$ as $x_{\tau}$ and $x_{\bar \tau}$ where $x_{\tau}$ consists of
the variables indexed by $\tau$ and $x_{\bar \tau}$ the variables
indexed by the complementary set $\bar \tau$. Similarly, the matrix
$A$ is partitioned as $A = [A_\tau, A_{\bar \tau}]$ and the cost
vector as $c = (c_{\tau}, c_{\bar \tau})$.  If $\sigma$ is a maximal
face of $\Delta_c$, then $A_{\sigma}$ is nonsingular and $Ax = b$ can
be written as $x_{\sigma} = A_{\sigma}^{-1}(b - A_{\bar \sigma}x_{\bar
  \sigma})$. Then $c \cdot x = c_{\sigma}(A_{\sigma}^{-1}(b - A_{\bar
  \sigma}x_{\bar \sigma})) + c_{\bar \sigma}x_{\bar \sigma} =
c_{\sigma}A_{\sigma}^{-1}b + (c_{\bar \sigma} -
c_{\sigma}A_{\sigma}^{-1}A_{\bar \sigma})x_{\bar \sigma}$. Let
$\tilde{c}_{{\bar \sigma}} := c_{\bar \sigma} -
c_{\sigma}A_{\sigma}^{-1}A_{\bar \sigma}$ and, for any face $\tau$ of
$\sigma$, let $\tilde{c}_{{\bar \tau}}$ be the extension of
$\tilde{c}_{{\bar \sigma}}$ to a vector in $\mathbb R^{|\bar \tau|}$
by adding zeros.

We now define a group relaxation of $IP_{A,c}(b)$ with respect to each 
face $\tau$ of $\Delta_c$.

\begin{definition} \label{grouprels}
The group relaxation of the integer program $IP_{A,c}(b)$ with respect
to the face $\tau$ of $\Delta_c$ is the program:
$$ G^{\tau}(b) = minimize \, \, \{ \tilde{c}_{{\bar \tau}} \cdot x_{\bar
\tau} \,\, : A_{\tau}x_{\tau} + A_{\bar \tau}x_{\bar \tau} = b, \,\,
x_{\bar \tau} \geq 0, (x_{\tau}, x_{\bar \tau}) \in \mathbb Z^n \}.$$
\end{definition}

Equivalently, $G^{\tau}(b) = minimize \, \, \{ \tilde{c}_{{\bar \tau}}
\cdot x_{\bar \tau} \,\, :\,\, A_{\bar \tau}x_{\bar \tau}\equiv b
\,\,(mod \,\, \mathbb Z A_{\tau}), \,\, x_{\bar \tau} \geq 0,
\,integer \}$ where $\mathbb Z A_{\tau}$ is the lattice generated by
the columns of $A_{\tau}$. Suppose $x_{\bar \tau}^{\ast}$ is an
optimal solution to the latter formulation. Since $\tau$ is a face of
$\Delta_c$, the columns of $A_{\tau}$ are linearly independent, and 
therefore the linear system $A_{\tau}x_{\tau} + A_{\bar
  \tau}x_{\bar \tau}^{\ast} = b$ has a unique solution. Solving this
system for $x_{\tau}$, the optimal solution $x_{\bar \tau}^{\ast}$ of
$G^{\tau}(b)$ can be uniquely lifted to the solution
$(x_{\tau}^{\ast}, x_{\bar \tau}^{\ast})$ of $Ax = b$. The formulation
of $G^{\tau}(b)$ in Definition~\ref{grouprels} shows that
$x_{\tau}^{\ast}$ is an integer vector. The group relaxation
$G^{\tau}(b)$ solves $IP_{A,c}(b)$ if and only if $x_{\tau}^{\ast}$ is
also non-negative.

The group relaxations of $IP_{A,c}(b)$ from Definition~\ref{grouprels}
contain among them the classical group relaxations of $IP_{A,c}(b)$
found in the literature. The program $G^{\sigma}(b)$, where
$A_{\sigma}$ is the optimal basis of the linear relaxation
$LP_{A,c}(b)$, is precisely Gomory's group relaxation of $IP_{A,c}(b)$
\cite{Gom65}.  The set of relaxations $G^{\tau}(b)$ as $\tau$ varies
among the subsets of this $\sigma$ are the {\em extended} group
relaxations of $IP_{A,c}(b)$ defined by Wolsey \cite{Wol71}.  Since
$\emptyset \in \Delta_c$, $G^{\emptyset}(b) = IP_{A,c}(b)$ is a group
relaxation of $IP_{A,c}(b)$, and hence $IP_{A,c}(b)$ will certainly be
solved by one of its extended group relaxations.  However, it is
possible to construct examples where a group relaxation $G^\tau(b)$
solves $IP_{A,c}(b)$, but $G^\tau(b)$ is neither Gomory's group
relaxation of $IP_{A,c}(b)$ nor one of its nontrivial extended Wolsey
relaxations (see Example \ref{long-chain}). Thus,
Definition~\ref{grouprels} typically creates more group relaxations
for each program in $IP_{A,c}$ than in the classical situation. This
has the obvious advantage that it increases the chance that
$IP_{A,c}(b)$ will be solved by some non-trivial relaxation, although
one may have to keep track of many more relaxations for each
program.  In Theorem~\ref{bounded}, we will prove that
Definition~\ref{grouprels} is the best possible in the sense that the
relaxations of $IP_{A,c}(b)$ defined there are precisely all the
bounded group relaxations of the program.

The goal in the rest of this section is to describe a useful
reformulation of the group problem $G^{\tau}(b)$ which is needed in
the rest of the paper and in the proof of Theorem~\ref{bounded}.
Given a sublattice $\Lambda$ of $\mathbb Z^n$, a cost vector $w \in
\mathbb R^n$ and a vector $v \in \mathbb N^n$, the {\em lattice
  program} defined by this data is
$$ minimize \,\,\{ w \cdot x \,\, : \,\, x \equiv v \,\, (mod \,\,
\Lambda),\,\, x \in \mathbb N^n \}.$$

Let $\mathcal L$ denote the $(n-d)$-dimensional saturated lattice $\{
x \in \mathbb Z^n : Ax = 0 \} \subseteq \mathbb Z^n$ and $u$ be a
feasible solution of the integer program $IP_{A,c}(b)$. Since
$IP_{A,c}(b)\,\, = \,\, minimize\,\, \{ c \cdot x \,\,: \,\, Ax = b
\,\,(= Au), \,\, x \in \mathbb N^n \}$ can be rewritten as
$minimize\,\,\{ c \cdot x \,\,: \,\,x -u \in \mathcal L, \,\, x \in
\mathbb N^n \}$, $IP_{A,c}(b)$ is equivalent to the lattice program
$$ minimize \,\, \{ c \cdot x \,\,: \,\,x \,\equiv \,u \,(mod \,\,
\mathcal L), \,\, x \in \mathbb N^n \}.$$

For $\tau \subseteq \{1, \ldots, n\}$, let $\pi_{\tau}$ be the
projection map from $\mathbb R^n \rightarrow \mathbb R^{|\bar \tau|}$
that kills all coordinates indexed by $\tau$. Then $\mathcal L_{\tau}
:= \pi_{\tau}(\mathcal L)$ is a sublattice of $\mathbb Z^{|\bar
  \tau|}$ that is isomorphic to $\mathcal L$: Clearly, $\pi_{\tau} :
\mathcal L \rightarrow \mathcal L_{\tau}$ is a surjection. If
$\pi_{\tau}(v) = \pi_{\tau}(v')$ for $v,v' \in \mathcal L$, then
$A_{\tau}v_{\tau} + A_{\bar \tau}v_{\bar \tau} = 0 = A_{\tau}v_{\tau}'
+ A_{\bar \tau}v_{\bar \tau}'$, implies that $A_{\tau}(v_{\tau} -
v_{\tau}') = 0$. Then $v_{\tau} = v_{\tau}'$ since the columns of
$A_{\tau}$ are linearly independent. Using this fact,
$G^{\tau}(b)$ can also be reformulated as a lattice program:\\

$\begin{array}{lll} G^{\tau}(b) & = & 
minimize \, \, \{ \tilde{c}_{{\bar \tau}} \cdot x_{\bar
\tau} \,\, : A_{\tau}x_{\tau} + A_{\bar \tau}x_{\bar \tau} = b, \,\,
x_{\bar \tau} \geq 0, (x_{\tau}, x_{\bar \tau}) \in \mathbb Z^n \}\\
& = & 
minimize\,\, \{\tilde c_{\bar \tau} \cdot x_{\bar \tau}  \,\,: \,\,
(x_{\tau},x_{\bar \tau})^t - (u_{\tau},u_{\bar \tau})^t \in \mathcal L, \,\,
x_{\bar \tau} \in \mathbb N^{|\bar \tau|} \} \\  
& = & 
minimize\,\, \{\tilde c_{\bar \tau} \cdot x_{\bar \tau}  \,\,: \,\,
x_{\bar \tau} - u_{\bar \tau} \in \mathcal L_{\tau}, \,\,
x_{\bar \tau} \in \mathbb N^{|\bar \tau|} \}\\
& = & 
minimize \,\, \{\tilde c_{\bar \tau} \cdot x_{\bar \tau} \,\,: \,\, 
x_{\bar \tau} \,\equiv \, \pi_{\tau}(u) \,(mod \,\, \mathcal
L_{\tau}), \,\, x_{\bar \tau} \in \mathbb N^{|\bar \tau|} \}. 
\end{array}$\\ 

Lattice programs were shown to be solved by Gr\"obner bases in
\cite{SWZ}. Theorem 5.3 in \cite{SWZ} gives a geometric
interpretation of these Gr\"obner bases in terms of corner polyhedra.
This paper was the first to make a connection between the theory of
group relaxations and commutative algebra (see \cite[\S 6]{SWZ}).
Special results are possible when the sublattice $\Lambda$ is of
finite index. In particular, the associated Gr\"obner bases are easier
to compute.

Since the $(n-d)$-dimensional lattice $\mathcal L \subset \mathbb Z^n$
is isomorphic to $\mathcal L_{\sigma} \subset \mathbb Z^{|\bar
  \sigma|}$ for $\sigma \in \Delta_c$, $\mathcal L_{\sigma}$ is of
finite index if and only if $\sigma$ is a maximal face of $\Delta_c$.
Hence the group relaxations $G^{\sigma}(b)$ as $\sigma$ varies over
the maximal faces of $\Delta_c$ are the easiest to solve among all
group relaxations of $IP_{A,c}(b)$. They contain among them Gomory's
group relaxation of $IP_{A,c}$. We give them a collective name in the
following definition.

\begin{definition} \label{gomory-relaxations}
The group relaxations $G^{\sigma}(b)$ of $IP_{A,c}(b)$, as $\sigma$
varies among the maximal faces of $\Delta_c$, are called the {\em
Gomory relaxations} of $IP_{A,c}(b)$.
\end{definition}

It is useful to reformulate $G^{\tau}(b)$ once again as follows. Let
$B \in \mathbb Z^{n \times (n-d)}$ be any matrix such that the columns
of $B$ generate the lattice $\mathcal L$, and let $u$ be a feasible
solution of $IP_{A,c}(b)$ as before. Then \\

$\begin{array}{lll} IP_{A,c}(b) 
  & = & minimize\,\, \{ c \cdot x \,\,:
  \,\,x \,- \,u \in \mathcal L, \,\, x \in \mathbb N^n \} \\
  & = & minimize \,\, \{ c \cdot x \,\,: \,\,x= u - Bz, \,\, x \geq 0,
  \,\, z \in
  \mathbb Z^{n-d} \}. \end{array}$\\
The last problem is equivalent to $minimize\,\, \{ c \cdot (u -
Bz)\,\,: \,\, Bz \leq u, \,\, z \in \mathbb Z^{n-d} \}$ and, therefore
$IP_{A,c}(b)$ is equivalent to the problem
\begin{equation} \label{ip}
minimize\,\, \{ (-cB) \cdot z \,\,: \,\, Bz \leq u, \,\, z \in \mathbb
Z^{n-d} \}. 
\end{equation}
There is a bijection between the set of feasible solutions of
(\ref{ip}) and the set of feasible solutions of $IP_{A,c}(b)$ via the 
isomorphism $z \mapsto u - Bz$. In particular, $0 \in \mathbb
R^{n-d}$ is feasible for (\ref{ip}) and it is the pre-image of $u$
under this map.

If $B^{\bar \tau}$ denotes the $|\bar \tau| \times (n-d)$
submatrix of $B$ obtained by deleting the rows indexed by $\tau$, then
$\mathcal L_{\tau} = \pi_{\tau}(\mathcal L) =
\{B^{\bar \tau}z \,\, : \,\, z \in \mathbb Z^{n-d} \}$.  
Using the same techniques as above, $G^{\tau}(b)$ can be reformulated
as  
$$
minimize\,\, \{ (-\tilde{c}_{{\bar \tau}}B^{\bar \tau}) \cdot
z \,\,:\,\, B^{\bar \tau}z \leq \pi_{\tau}(u), \,\, z \in \mathbb
Z^{n-d} \}.
$$

Since $\tilde{c}_{{\bar \tau}} = \pi_{\tau}(c -
c_{\sigma}A_{\sigma}^{-1}A)$ for any maximal face $\sigma$ of
$\Delta_c$ containing $\tau$ and the support of $c -
c_{\sigma}A_{\sigma}^{-1}A$ is contained in $\bar \tau$, 
$\tilde{c}_{{\bar \tau}}B^{\bar \tau} = (c -
c_{\sigma}A_{\sigma}^{-1}A)B = cB$ since $AB = 0$. Hence $G^{\tau}(b)$
is equivalent to 
\begin{equation} \label{gp}
minimize\,\, \{ (-cB) \cdot z \,\,:\,\, B^{\bar \tau}z \leq
\pi_{\tau}(u), \,\, z \in \mathbb Z^{n-d} \}.
\end{equation}
The feasible solutions to (\ref{ip}) are the lattice points in the
rational polyhedron $P_{u} := \{ z \in \mathbb R^{n-d} : Bz \leq u
\}$, and the feasible solutions to (\ref{gp}) are the lattice points
in the relaxation $P_{u}^{\bar \tau} := \{ z \in \mathbb R^{n-d} :
B^{\bar \tau} z \leq \pi_{\tau}(u) \}$ of $P_{u}$ obtained by deleting
the inequalities indexed by $\tau$. In theory, one could define group
relaxations of $IP_{A,c}(b)$ with respect to any $\tau \subseteq \{1,
\ldots, n\}$.  The following theorem illustrates the completeness of
Definition~\ref{grouprels}.

\begin{theorem} \label{bounded}
The group relaxation $G^{\tau}(b)$ of $IP_{A,c}(b)$ 
has a finite optimal solution if and  
only if $\tau \subseteq \{1, \ldots, n\}$ is a face of $\Delta_c$.
\end{theorem}

\begin{proof}
Since all data are integral it suffices to prove that the linear
relaxation $$minimize\,\, \{ (-cB) \cdot z \,\, : \,\, z \in P_{u}^{\bar
\tau} \}$$ is bounded if and only if $\tau \in \Delta_c$.  

If $\tau$ is a face of $\Delta_c$ then there exists $y \in \mathbb
R^d$ such that $yA_{\tau} = c_{\tau}$ and $yA_{\bar \tau} < c_{\bar
\tau}$. Using the fact that $A_{\tau}B^{\tau} + A_{\bar \tau}B^{\bar
\tau} = 0$ we see that $cB = c_{\tau}B^{\tau} + c_{\bar \tau}B^{\bar
\tau} = yA_{\tau}B^{\tau} + c_{\bar \tau}B^{\bar \tau} = y(-A_{\bar
\tau}B^{\bar \tau}) + c_{\bar \tau}B^{\bar \tau} = (c_{\bar \tau} -
yA_{\bar \tau})B^{\bar \tau}$. This implies that $cB$ is a positive
linear combination of the rows of $B^{\bar \tau}$ since $c_{\bar \tau}
- yA_{\bar \tau} > 0$. Hence $cB$ lies in the polar of $\{z \in
\mathbb R^{n-d} : B^{\bar \tau}z \leq 0 \}$ which is the recession
cone of $P_{u}^{\bar \tau}$ proving that the linear program $minimize 
\,\, \{ (-cB) \cdot z \,\, : \,\, z \in P_{u}^{\bar \tau} \}$ is bounded.

The linear program $minimize\,\, \{ (-cB) \cdot z \,\, : \,\, z \in
P_{u}^{\bar \tau} \}$ is feasible since $0$ is a feasible solution. If
it is bounded as well then $minimize\,\, \{ c_{\tau}x_{\tau} + c_{\bar
  \tau}x_{\bar \tau} \,\, : \,\, A_{\tau}x_{\tau} + A_{\bar
  \tau}x_{\bar \tau} = b, \, x_{\bar \tau} \geq 0 \}$ is feasible and
bounded. As a result, the dual of the latter program $maximize\,\, \{
y \cdot b \,\, : \,\, yA_{\tau} = c_{\tau}, \, yA_{\bar \tau} \leq
c_{\bar \tau} \}$ is feasible. This shows that a superset of $\tau$ is
a face of $\Delta_c$ which implies that $\tau \in \Delta_c$ since
$\Delta_c$ is a triangulation.
\end{proof}

\section{Associated Sets}
The group relaxation $G^{\tau}(b)$ (seen as (\ref{gp})) solves the
integer program $IP_{A,c}(b)$ (seen as (\ref{ip})) if and only if both
programs have the same optimal solution $z^{\ast} \in \mathbb
Z^{n-d}$. If $G^{\tau}(b)$ solves $IP_{A,c}(b)$ then $G^{\tau'}(b)$
also solves $IP_{A,c}(b)$ for every $\tau' \subset \tau$ since
$G^{\tau'}(b)$ is a {\em stricter} relaxation of $IP_{A,c}(b)$ than
$G^{\tau}(b)$. For the same reason, one would expect that
$G^{\tau}(b)$ is easier to solve than $G^{\tau'}(b)$.  Therefore, the
most useful group relaxations of $IP_{A,c}(b)$ are those indexed by
the {\em maximal} elements in the subcomplex of $\Delta_c$ consisting
of all faces $\tau$ such that $G^{\tau}(b)$ solves $IP_{A,c}(b)$. The
following definition isolates such relaxations.

\begin{definition} \label{associatedset}
  A face $\tau$ of the regular triangulation $\Delta_c$ is
  an {\em associated set} of $IP_{A,c}$ (or is associated to $IP_{A,c}$)
  if for some $b \in \mathbb N A$, $G^{\tau}(b)$ solves $IP_{A,c}(b)$
  but $G^{\tau'}(b)$ does not for all faces $\tau'$ of $\Delta_c$ such
  that $\tau \subset \tau'$.
\end{definition}

The associated sets of $IP_{A,c}$ carry all the information about all
the group relaxations needed to solve the programs in $IP_{A,c}$.  In
this section we will develop tools to understand these sets. We start
by considering the set $\mathcal O_c \subset \mathbb N^n$ of all the
optimal solutions of all programs in $IP_{A,c}$. A basic result
in the algebraic study of integer programming is that $\mathcal O_c$ is
an {\em order ideal} or {\em down set} in $\mathbb N^n$, i.e., if $u
\in \mathcal O_c$ and $v \leq u, \,\, v \in \mathbb N^n$, then $v \in
\mathcal O_c$. One way to prove this is to show that the complement
$\mathcal N_c := \mathbb N^n \backslash \mathcal O_c$ has the property
that if $v \in \mathcal N_c$ then $v + \mathbb N^n \subseteq \mathcal
N_c$. Every lattice point in $\mathbb N^n$ is a feasible solution to
a unique program in $IP_{A,c}$ ($u \in \mathbb N^n$ is feasible for
$IP_{A,c}(Au)$). Hence, $\mathcal N_c$ is the set of all non-optimal
solutions of all programs in $IP_{A,c}$. A set $P \subseteq \mathbb
N^n$ with the property that $p + \mathbb N^n \subseteq P$ whenever $p
\in P$ has a finite set of minimal elements. Hence there exists
$\alpha_1, \ldots, \alpha_t \in \mathcal N_c$ such that
$$
\mathcal N_c = \bigcup_{i=1}^{t} (\alpha_i + \mathbb N^n).$$
As a
result, $\mathcal O_c$ is completely specified by the finitely many
``generators'' $\alpha_1,\ldots, \alpha_t$ of its complement $\mathcal
N_c$. See \cite{Tho} for proofs of these assertions.

\begin{example} \label{knapsack} 
Consider the family of knapsack problems:
$$minimize \,\,\{10000x_1 + 100x_2 + x_3 \,\,:\,\, 2x_1 + 5x_2 + 8x_3
= b,\,\, (x_1,x_2,x_3) \in \mathbb N^3 \}$$
as $b$ varies in the semigroup $\mathbb N [2 \,\,5\,\,8]$. The set
$\mathcal N_c$ is generated by the vectors 
$$(0,8,0),(1,0,1),(1,6,0),(2,4,0),(3,2,0),\,\,{\text and}\,\,(4,0,0)$$
which means that $\mathcal N_c = ((0,8,0) + \mathbb N^3) \cup \cdots
\cup ((4,0,0) + \mathbb N^3).$ Figure~\ref{f4} is a picture of
$\mathcal N_c$ (created by Ezra Miller). The white points are its
generators.  One can see that $\mathcal O_c$ consists of finitely many
points of the form $(p,q,0)$ where $p \geq 1$ and the eight ``lattice
lines'' of points $(0,i,\ast)$, $i = 0, \ldots, 7$.
\begin{figure}
\epsfig{file=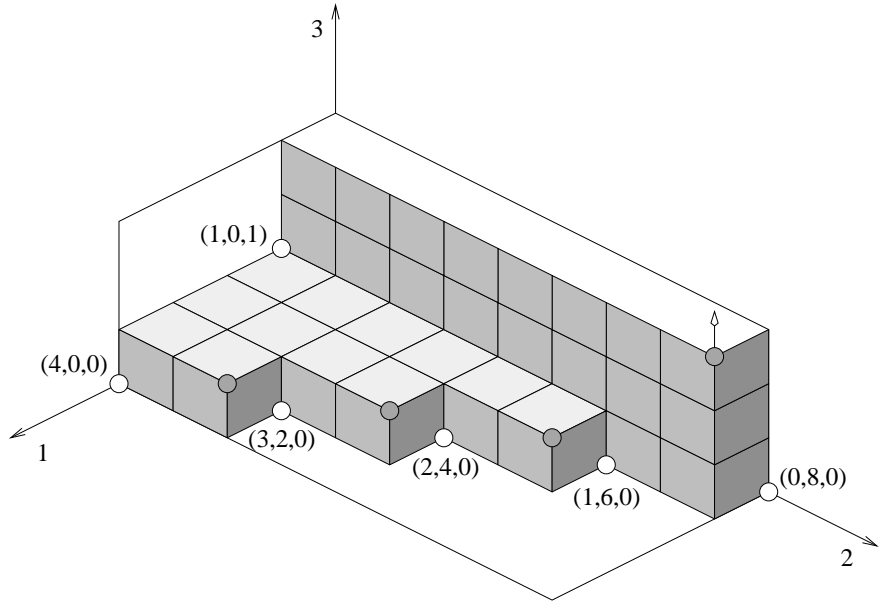, height=8cm}
\caption{The set of non-optimal solutions $\mathcal N_c$ for
  Example~\ref{knapsack}.} 
\label{f4}
\end{figure}
\qed
\end{example}

For the purpose of computations, it is most effective to think of
$\mathcal N_c$ and $\mathcal O_c$ algebraically. A {\em monomial}
$x^u$ in the {\em polynomial ring } $S := k[x_1, \ldots, x_n]$ is a
product $x^u = x_1^{u_1}x_2^{u_2}\cdots x_n^{u_n}$ where $u = (u_1,
\ldots, u_n) \in \mathbb N^n$. We assume that $k$ is a field, say the
set of rational numbers. For a scalar $k_u \in k$ and a monomial $x^u$
in $S$, we call $k_ux^u$ a {\em term} of $S$. A {\em polynomial} $f =
\sum k_ux^u$ in $S$ is a combination of finitely many terms in $S$.  A
subset $I$ of $S$ is an {\em ideal} of $S$ if (1) $I$ is closed under
addition, i.e., $f,g \in I \Rightarrow f+g \in I$ and (2) if $f \in I$
and $g \in S$ then $fg \in I$. We say that $I$ is generated by the
polynomials $f_1, \ldots, f_t$, denoted as $I = \langle f_1, \ldots,
f_t \rangle$, if $I = \{\sum_{i=1}^{t} f_ig_i \,\, : \,\, g_i \in S
\}$.  By Hilbert's basis theorem, every ideal in $S$ has a finite
generating set. An ideal $M$ in $S$ is called a {\em monomial ideal}
if it is generated by monomials, i.e., $M = \langle x^{v_1}, \ldots,
x^{v_t} \rangle$ for monomials $x^{v_1}, \ldots, x^{v_t}$ in $S$. The
monomials that do not lie in $M$ are called the {\em standard
  monomials} of $M$.  The {\em cost} of a term $k_ux^u$ with respect
to a vector $c \in \mathbb R^n$ is the dot product $c \cdot u$. The
{\em initial term} of a polynomial $f = \sum k_u x^u \in S$ with
respect to $c$, denoted as $in_c(f)$, is the sum of all terms in $f$
of maximal cost. For any ideal $I \subset S$, the {\em initial ideal}
of $I$ with respect to $c$, denoted as $in_c(I)$, is the ideal
generated by all the initial terms $in_c(f)$ of all polynomials $f$ in
$I$. These concepts come from the theory of {\em Gr\"obner bases} for
polynomial ideals. See \cite{CLO} for an introduction.

The {\em toric ideal} of the matrix $A$,
denoted as $I_A$, is the {\em binomial} ideal in $S$ defined as:
$$
I_A := \langle x^u - x^v : \,\, u, v \in \mathbb{N}^n \text{ and }
Au = Av \rangle. $$
Toric ideals provide the link between integer
programming and Gr\"obner basis theory. See \cite{Stu} and
\cite{Tho98} for an introduction to this area of research. This
connection yields the following basic facts that we state without
proofs. (Recall that the cost vector $c$ of $IP_{A,c}$ was assumed to
be generic in the sense that each program in $IP_{A,c}$ has a unique
optimal solution.)

\begin{lemma} \label{grobnerfacts} \cite{Stu}
  (i) If $c$ is generic, then the initial ideal $in_c(I_A)$ is a
  monomial ideal.\\
  (ii) A lattice point $u$ is non-optimal for the integer program
  $IP_{A,c}(Au)$, or equivalently, $u \in \mathcal N_c$, if and only
  if $x^u$ lies in the initial ideal $in_c(I_A)$. In other words, a
  lattice point $u$ lies in $\mathcal O_c$ if and only if $x^u$ is a
  standard monomial of $in_c(I_A)$.\\
  (iii) The {\bf reduced Gr\"obner basis} $\mathcal G_c$ of $I_A$ with
  respect to $c$ is the unique minimal test set for the family of
  integer programs $IP_{A,c}$.\\
  (iv) If $u$ is a feasible solution of $IP_{A,c}(b)$, and
  $x^{u^{\ast}}$ is the unique {\bf normal form} of $x^u$ with respect
  to $\mathcal G_c$, then $u^{\ast}$ is the optimal solution of
  $IP_{A,c}(b)$.
\end{lemma}

We do not elaborate on parts (iii) and (iv) of
Lemma~\ref{grobnerfacts}. They are not needed for what follows and are
included for completeness. Since $c$ is generic,
Lemma~\ref{grobnerfacts} (ii) implies that there is a bijection
between the lattice points of $\mathcal O_c$ and the semigroup
$\mathbb N A$ via the map $\phi_A : \mathcal O_c \rightarrow \mathbb N
A$ such that $u \mapsto Au$. The inverse of $\phi_A$ sends a vector $b
\in \mathbb N A$ to the optimal solution of $IP_{A,c}(b)$.\\ 

\noindent {\bf Example~\ref{knapsack} continued.}
In this example, the toric ideal $I_A = \langle x_1^4-x_3,
x_2^2-x_1x_3\rangle$ and its initial ideal with respect to the cost
vector $c = (10000,100,1)$ is $$in_c(I_A) = \langle x_2^8, \,x_1x_3,
\,x_1x_2^6, \,x_1^2x_2^4, \,x_1^3x_2^2, \,x_1^4 \rangle.$$
Note that
the exponent vectors of the generators of $in_c(I_A)$ are 
the generators of $\mathcal N_c$. \qed\\

We will now describe a certain decomposition of the set $\mathcal O_c$
which in turn will shed light on the associated sets of $IP_{A,c}$.
For $u \in \mathbb N^n$, consider $Q_u := \{ z \in \mathbb R^{n-d}
\,\, : \,\, Bz \leq u, \, (-cB) \cdot z \leq 0 \}$ and its relaxation
$Q_u^{\bar \tau} := \{ z \in \mathbb R^{n-d} \,\, : \,\, B^{\bar \tau}
z \leq \pi_{\tau}(u), \, (-cB) \cdot z \leq 0 \}$ where $B, B^{\bar \tau}$
are as in (\ref{ip}) and (\ref{gp}) and $\tau \in \Delta_c$. By
Theorem~\ref{bounded}, both $Q_u$ and $Q_u^{\bar \tau}$ are polytopes.
Notice that if $\pi_{\tau}(u) = \pi_{\tau}(u')$ for two distinct
vectors $u,u' \in \mathbb N^n$, then $Q_u^{\bar \tau} = Q_{u'}^{\bar
  \tau}$.

\begin{lemma} \label{emptypolys}
(i) A lattice point $u$ is in  $\mathcal O_c$ if and only
if $Q_{u} \cap \mathbb Z^{n-d} = \{0\}$.\\
(ii) If $u \in \mathcal O_c$, then the group relaxation $G^{\tau}(Au)$
solves the integer program $IP_{A,c}(Au)$ if and only if $Q_{u}^{\bar
\tau} \cap \mathbb Z^{n-d} = \{0\}$.
\end{lemma}

\begin{proof} (i) The lattice point $u$ belongs to $\mathcal O_c$ if
and only if $u$ is the optimal solution to $IP_{A,c}(Au)$ which is
equivalent to $0 \in \mathbb Z^{n-d}$ being the optimal solution to
the reformulation (\ref{ip}) of $IP_{A,c}(Au)$. Since $c$ is generic,
the last statement is equivalent to $Q_{u} \cap \mathbb Z^{n-d} =
\{0\}$. The second statement follows from (i) and the fact that
(\ref{gp}) solves (\ref{ip}) if and only if they have the same optimal
solution. 
\end{proof}

In order to state the coming results, it is convenient to assume that
the vector $u$ in (\ref{ip}) and (\ref{gp}) is the optimal solution to
$IP_{A,c}(b)$. For an element $u \in \mathcal O_c$ and a face $\tau$
of $\Delta_c$ let $S(u,\tau)$ be the {\em affine} semigroup $u +
\mathbb N (e_i : i \in \tau) \subseteq \mathbb N^n$ where $e_i$
denotes the $i$th unit vector of $\mathbb R^n$. Note that $S(u,\tau)$
is not a semigroup if $u \neq 0$ (since $0 \not \in S(u,\tau)$), but
is a translation of the semigroup $\mathbb N (e_i : i \in \tau)$. We
use the adjective {\em affine} here as in an {\em affine subspace}
which is not a subspace but the translation of one.  Note that if $v
\in S(u,\tau)$, then $\pi_{\tau}(v) = \pi_{\tau}(u)$.

\begin{lemma} \label{stdpair-grouprel} 
For $u \in \mathcal O_c$ and a face $\tau$ of $\Delta_c$, the affine
semigroup $S(u,\tau)$ is contained in $\mathcal O_c$ if and only if
$G^{\tau}(Au)$ solves $IP_{A,c}(Au)$.
\end{lemma}

\begin{proof} 
Suppose $S(u,\tau) \subseteq \mathcal O_c$. Then by
Lemma~\ref{emptypolys} (i), for all $v \in S(u,\tau)$, $$Q_v = \{ z
\in \mathbb R^{n-d} \,\, : \,\, B^{\tau}z \leq \pi_{\bar \tau}(v), \,
B^{\bar \tau}z \leq \pi_{\tau}(u), \, (-cB) \cdot z \leq 0 \} \cap \mathbb
Z^{n-d} = \{0\}.$$ Since $\pi_{\bar \tau}(v)$ can be any vector in
$\mathbb N^{|\tau|}$, $Q_u^{\bar \tau} \cap \mathbb Z^{n-d} =
\{0\}$. Hence, by Lemma~\ref{emptypolys} (ii), $G^{\tau}(Au)$ solves
$IP_{A,c}(Au)$.

If $v \in S(u,\tau)$, then $\pi_{\tau}(u) =
\pi_{\tau}(v)$, and hence $Q_u^{\bar \tau} = Q_v^{\bar
\tau}$. Therefore, if $G^{\tau}(Au)$ solves $IP_{A,c}(Au)$, then
$\{0\} = Q_u^{\bar \tau} \cap \mathbb Z^{n-d} = Q_v^{\bar \tau} \cap
\mathbb Z^{n-d}$ for all $v \in S(u,\tau)$.  Since $Q_v^{\bar \tau}$
is a relaxation of $Q_v$, $Q_v \cap \mathbb Z^{n-d} = \{0\}$ for all
$v \in S(u,\tau)$ and hence by Lemma~\ref{emptypolys} (i), $S(u,\tau) 
\subseteq \mathcal O_c$.
\end{proof}

\begin{lemma} \label{solvesall}
For $u \in \mathcal O_c$ and a face $\tau$ of $\Delta_c$,
$G^{\tau}(Au)$ solves $IP_{A,c}(Au)$ if and only if $G^{\tau}(Av)$
solves $IP_{A,c}(Av)$ for all $v \in S(u,\tau)$.
\end{lemma}

\begin{proof}
  If $v \in S(u,\tau)$ and $G^{\tau}(Au)$ solves $IP_{A,c}(Au)$, then
  as seen before, $\{0\} = Q_u^{\bar \tau} \cap \mathbb
  Z^{n-d} = Q_v^{\bar \tau} \cap \mathbb Z^{n-d}$ for all $v \in
  S(u,\tau)$. By Lemma~\ref{emptypolys} (ii), $G^{\tau}(Av)$ solves
  $IP_{A,c}(Av)$ for all $v \in S(u,\tau)$. The converse holds for the
  trivial reason that $u \in S(u,\tau)$.
\end{proof}

\begin{corollary}
For $u \in \mathcal O_c$ and a face $\tau$ of $\Delta_c$, the affine
semigroup $S(u,\tau)$ is contained in $\mathcal O_c$ if and only if
$G^{\tau}(Av)$ solves $IP_{A,c}(Av)$ for all $v \in S(u,\tau)$.
\end{corollary}

Since $\pi_{\tau}(u)$ determines the polytope $Q_u^{\bar \tau} =
Q_v^{\bar \tau}$ for all $v \in S(u,\tau)$, we could have assumed that 
$supp(u) \subseteq {\bar \tau}$ in Lemmas~\ref{stdpair-grouprel} and
\ref{solvesall}.

\begin{definition} 
For $\tau \in \Delta_c$ and $u \in \mathcal O_c$, $(u,\tau)$
is called an {\em admissible pair} of $\mathcal O_c$ if \\ 
\indent (i) the support of $u$ is contained in $\bar \tau$, and \\
\indent (ii) $S(u,\tau) \subseteq \mathcal O_c$ or equivalently,
$G^{\tau}(Av)$ solves $IP_{A,c}(Av)$ for all $v \in S(u,\tau)$.\\ 
An admissible pair $(u, \tau)$ is a {\em standard pair} of $\mathcal
O_c$ if the affine semigroup $S(u,\tau)$ is not properly contained
in $S(v,\tau')$ where $(v,\tau')$ is another admissible pair of
$\mathcal O_c$.\\ 
\end{definition}

\noindent{\bf Example~\ref{knapsack} continued.} 
From Figure~\ref{f4}, one can see that the standard pairs of 
$\mathcal O_c$ are 
\begin{center}
$\begin{array}{ll}
((1,0,0),\emptyset) & ((1,3,0),\emptyset)\\
((2,0,0),\emptyset) & ((2,3,0),\emptyset)\\
((3,0,0),\emptyset) & ((1,4,0),\emptyset)\\
((1,1,0),\emptyset) & ((1,5,0),\emptyset)\\
((2,1,0),\emptyset) &\\
((3,1,0),\emptyset) &\\
((1,2,0),\emptyset) &\\
((2,2,0),\emptyset) &\\
\end{array}$
\hspace{.5in} and 
\hspace{.5in}
$\begin{array}{l}
((0,0,0),\{3\})\\
((0,1,0),\{3\})\\
((0,2,0),\{3\})\\ 
((0,3,0),\{3\})\\
((0,4,0),\{3\})\\
((0,5,0),\{3\})\\
((0,6,0),\{3\})\\
((0,7,0),\{3\})\\
\end{array}$
\end{center}
\qed

\begin{definition}\label{stdpolys}
For a face $\tau$ of $\Delta_c$ and a lattice point $u \in \mathbb
N^n$, we say that the polytope $Q_{u}^{\bar \tau}$ is a {\em standard
polytope} of $IP_{A,c}$ if $Q_{u}^{\bar \tau} \cap \mathbb Z^{n-d} =
\{0\}$ and every relaxation of $Q_{u}^{\bar \tau}$ obtained by
removing an inequality in $B^{\bar \tau}z \leq \pi_{\tau}(u)$ contains
a non-zero lattice point.
\end{definition}

\begin{figure}
\epsfig{file=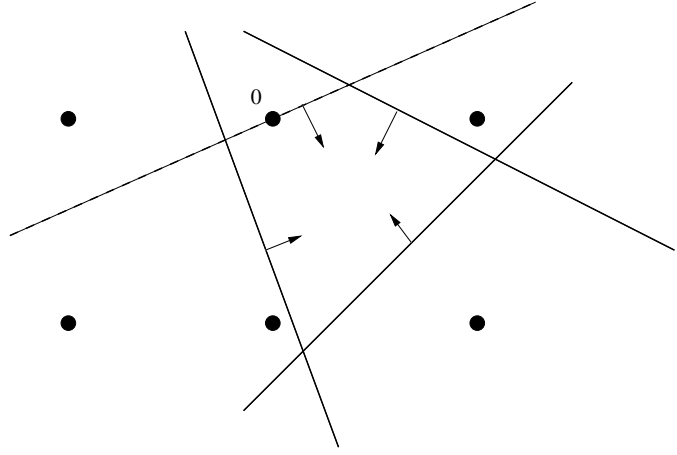, height=6cm}
\caption{A standard polytope.}
\label{f5}
\end{figure}

Figure~\ref{f5} is a diagram of a standard polytope $Q_{u}^{\bar
  \tau}$. The dashed line is the boundary of the half space $(-cB)
\cdot z \leq 0$ while the other lines are the boundaries of the
halfspaces given by the inequalities in $B^{\bar \tau}z \leq
\pi_{\tau}(u)$. The origin is the only lattice point in the polytope
and if any inequality in $B^{\bar \tau}z \leq \pi_{\tau}(u)$ is
removed, a lattice point will enter the relaxation.

We re-emphasize that if $Q_{u}^{\bar \tau}$ is a standard polytope,
then $Q_{u'}^{\bar \tau}$ is the same standard polytope if
$\pi_{\tau}(u) = \pi_{\tau}(u')$. Hence the same standard polytope can
be indexed by infinitely many $u \in \mathbb N^n$. We now state the
main result of this section which characterizes associated sets in
terms of standard pairs and standard polytopes.

\begin{theorem} \label{stdpairs-equivs}
The following statements are equivalent:\\ 
(i) The admissible pair $(u,\tau)$ is a standard pair of $\mathcal
O_c$.\\ 
(ii) The polytope $Q_{u}^{\bar \tau}$ is a standard polytope of
$IP_{A,c}$.\\ 
(iii) The face $\tau$ of $\Delta_c$ is associated to $IP_{A,c}$.
\end{theorem}

\begin{proof} $(i) \Leftrightarrow (ii)$: The admissible pair $(u,
\tau)$ is standard if and only if for every $i \in \bar \tau$, there
exists some positive integer $m_i$ and a vector $v \in S(u,\tau)$ such
that $v + m_ie_i \in \mathcal N_c$. (If this condition did not hold
for some $i \in \bar \tau$, then $(u', \tau \cup \{i\})$ would be an
admissible pair of $\mathcal O_c$ such that $S(u',\tau \cup \{i\})$
contains $S(u,\tau)$ where $u'$ is obtained from $u$ by setting the
$i$th component of $u$ to zero. Conversely, if the condition holds for
an admissible pair then the pair is standard.) Equivalently, for each
$i \in \bar \tau$, there exists a positive integer $m_i$ and a $v \in
S(u,\tau)$ such that $Q_{v+m_ie_i}^{\bar \tau} = Q_{u + m_ie_i}^{\bar
\tau}$ contains at least two lattice points. In other words, the
removal of the inequality indexed by $i$ from the inequalities in
$B^{\bar \tau}z \leq \pi_{\tau}(u)$ will bring an extra lattice point
into the corresponding relaxation of $Q_u^{\bar \tau}$. This is
equivalent to saying that $Q_u^{\bar \tau}$ is a standard polytope of
$IP_{A,c}$.\\ 

\noindent $(i) \Leftrightarrow (iii)$: Suppose $(u,\tau)$ is a
standard pair of $\mathcal O_c$. Then $S(u,\tau) \subseteq \mathcal
O_c$ and $G^{\tau}(Au)$ solves $IP_{A,c}(Au)$ by
Lemma~\ref{stdpair-grouprel}. Suppose $G^{\tau'}(Au)$ solves
$IP_{A,c}(Au)$ for some face $\tau' \in \Delta_c$ such that $\tau
\subset \tau'$. Lemma~\ref{stdpair-grouprel} then implies that
$S(u,\tau')$ lies in $\mathcal O_c$. This contradicts the fact that
$(u,\tau)$ was a standard pair of $\mathcal O_c$ since $S(u,\tau)$ is
properly contained in $S(\hat u,\tau')$ corresponding to the
admissible pair $(\hat u, \tau')$ where $\hat u$ is obtained from $u$
by setting $u_i = 0$ for all $i \in \tau' \backslash \tau$.

To prove the converse, suppose $\tau$ is associated to $IP_{A,c}$.
Then there exists some $b \in \mathbb N A$ such that $G^{\tau}(b)$
solves $IP_{A,c}(b)$ but $G^{\tau'}(b)$ does not for all faces $\tau'$
of $\Delta_c$ containing $\tau$. Let $u$ be the unique optimal
solution of $IP_{A,c}(b)$. By Lemma~\ref{stdpair-grouprel}, $S(u,\tau)
\subseteq \mathcal O_c$.  Let $\hat u \in \mathbb N^n$ be obtained
from $u$ by setting $u_i = 0$ for all $i \in \tau$. Then $G^{\tau}(A
\hat u)$ solves $IP_{A,c}(A \hat u)$ since $Q_u^{\bar \tau} = Q_{\hat
  u}^{\bar \tau}$. Hence $S(\hat u, \tau) \subseteq \mathcal O_c$ and
$(\hat u, \tau)$ is an admissible pair of $\mathcal O_c$. Suppose
there exists another admissible pair $(w,\sigma)$ such that $S(\hat u,
\tau) \subset S(w,\sigma)$. Then $\tau \subseteq \sigma$. If $\tau =
\sigma$ then $S(\hat u, \tau)$ and $S(w,\sigma)$ are both orthogonal
translates of $\mathbb N (e_i \,\,: \,\,i \in \tau )$ and hence
$S(\hat u, \tau)$ cannot be properly contained in $S(w,\sigma)$.
Therefore, $\tau$ is a proper subset of $\sigma$ which implies that
$S(\hat u, \sigma) \subseteq \mathcal O_c$. Then, by
Lemma~\ref{stdpair-grouprel}, $G^{\sigma}(A \hat u)$ solves
$IP_{A,c}(A \hat u)$ which contradicts that $\tau$ was an associated
set of $IP_{A,c}$.
\end{proof}

\noindent{\bf Example~\ref{knapsack} continued.} In
Example~\ref{knapsack} we can choose $B$ to be the $3 \times 2$ matrix
  $$B = \left[ \begin{array}{cc} -1 & 4 \\ 2 & 0 \\ -1
  & -1 \end{array} \right ].$$
The standard polytope defined by the standard pair
$((1,0,0),\emptyset)$ is hence 
$$ \{ (z_1,z_2) \in \mathbb R^2 \,\,:\,\, -z_1 + 4z_2 \leq 1,\,\, 
2z_1 \leq 0,\,\,  -z_1-z_2 \leq 0,\,\, 9801z_1 - 40001z_2 \leq 0 \}$$
while the standard polytope defined by the standard pair
$((0,2,0),\{3\})$ is 
$$ \{ (z_1,z_2) \in \mathbb R^2 \,\,:\,\, -z_1 + 4z_2 \leq 0,\,\, 
2z_1 \leq 2,\,\, 9801z_1 - 40001z_2 \leq 0 \}.$$
The associated sets of $IP_{A,c}$ in this example are $\emptyset$ and 
$\{3\}$. There are twelve quadrangular and eight triangular
standard polytopes for this family of knapsack problems. \qed\\

Standard polytopes were introduced in \cite{HT2}, and the equivalence
of parts (i) and (ii) of Theorem~\ref{stdpairs-equivs} was proved in
\cite[Theorem 2.5]{HT2}. Under the linear map $\phi_A : \mathbb N^n
\rightarrow \mathbb N A$ where $u \mapsto Au$, the affine
semigroup $S(u,\tau)$ where $(u,\tau)$ is a standard pair of $\mathcal
O_c$ maps to the affine semigroup $Au + \mathbb N A_{\tau}$ in
$\mathbb N A$. Since every integer program in $IP_{A,c}$ is solved by
one of its group relaxations, $\mathcal O_c$ is covered by the affine
semigroups corresponding to its standard pairs.  We call this cover
and its image in $\mathbb N A$ under $\phi_A$ the {\em standard pair
  decompositions} of $\mathcal O_c$ and $\mathbb N A$, respectively.
Since standard pairs of $\mathcal O_c$ are determined by the standard
polytopes of $IP_{A,c}$, the standard pair decomposition of $\mathcal
O_c$ is unique. The terminology used above has its origins in
\cite{STV} which introduced the {\em standard pair decomposition} of a
{\em monomial ideal}. The specialization to integer programming appear
in \cite{HT2}, \cite{HT1} and \cite[\S 12.D]{Stu}. The following
theorem shows how the standard pair decomposition of $\mathcal O_c$
dictates which group relaxations solve which programs in $IP_{A,c}$.

\begin{theorem} \label{cover=solve}
Let $v$ be the optimal solution of the integer program
$IP_{A,c}(b)$. Then the group relaxation $G^{\tau}(Av)$ solves
$IP_{A,c}(Av)$ if and only if there is some standard pair $(u,\tau')$
of $\mathcal O_c$ with $\tau \subseteq \tau'$ such that $v$ belongs to
the affine semigroup $S(u,\tau')$.
\end{theorem}

\begin{proof}
Suppose $v$ lies in $S(u,\tau')$ corresponding to the standard pair
$(u,\tau')$ of $\mathcal O_c$. Then $S(v,\tau') \subseteq \mathcal
O_c$ which implies that $G^{\tau'}(Av)$ solves $IP_{A,c}(Av)$ by
Lemma~\ref{stdpair-grouprel}. Hence $G^{\tau}(Av)$ also solves
$IP_{A,c}(Av)$ for all $\tau \subseteq \tau'$.

To prove the converse, suppose $\tau'$ is a maximal element in the
subcomplex of all faces $\tau$ of $\Delta_c$ such that $G^{\tau}(Av)$
solves $IP_{A,c}(Av)$. Then $\tau'$ is an associated set of
$IP_{A,c}$. In the proof of $(iii) \Rightarrow (i)$ in
Theorem~\ref{stdpairs-equivs}, we showed that $(\hat v, \tau')$ is a
standard pair of $\mathcal O_c$ where $\hat v$ is obtained from $v$ by
setting $v_i = 0$ for all $i \in \tau'$. Then $v \in S(\hat v,
\tau')$.
\end{proof}

\noindent{\bf Example~\ref{knapsack} continued.} The eight standard
pairs of $\mathcal O_c$ of the form $(\ast, \{3\})$, map to the 
eight affine semigroups:
$$\mathbb N [8],\, (5 + \mathbb N [8]),\, (10 + \mathbb N [8]),\, (15
+ \mathbb N [8]),\, (20 + \mathbb N [8]),\, (25 + \mathbb N [8]),\,
(30 + \mathbb N [8]) \,\, {\text and}\,\,(35 + \mathbb N [8])$$
contained in $\mathbb N A = \mathbb N \,[2,5,8] \subset \mathbb N$.
For all right hand side vectors $b$ in the union of these sets, the
integer program $IP_{A,c}(b)$ can be solved by the group relaxation
$G^{\{3\}}(b)$.  The twelve standard pairs of the from
$(\ast,\emptyset)$ map to the remaining finitely many points
$$
2,4,6,7,9,11,12,14,17,19,22\,\,{\text and}\,\,27$$
of $\mathbb N
\,[2,5,8]$.  If $b$ is one of these points, then $IP_{A,c}(b)$ can
only be solved as the full integer program. In this example, the
regular triangulation $\Delta_c = \{\{3\}\}$. Hence $G^{\{3\}}(b)$ is
a Gomory relaxation of $IP_{A,c}(b)$. \qed \\

For most $b \in \mathbb N A$, the program $IP_{A,c}(b)$ is solved by
one of its Gomory relaxations, or equivalently, by
Theorem~\ref{cover=solve}, the optimal solution $v$ of $IP_{A,c}(b)$
lies in $S(\ast,\sigma)$ for some standard pair $(\ast,\sigma)$ where
$\sigma$ is a maximal face of $\Delta_c$. For mathematical versions of
this informal statement, see \cite[Proposition 12.16]{Stu} and
\cite[Theorems 1 and 2]{Gom65}. Roughly speaking, these right hand
sides are away from the boundary of $cone(A)$. (This was seen in
Example~\ref{knapsack} above, where for all but twelve right hand
sides, $IP_{A,c}(b)$ was solvable by the Gomory relaxation
$G^{\{3\}}(b)$.  Further, these right hand sides were toward the
boundary of $cone(A)$, the origin in this one-dimensional case.)  For
the remaining right hand sides, $IP_{A,c}(b)$ can only be solved by
$G^{\tau}(b)$ where $\tau$ is a lower dimensional face of $\Delta_c$ -
possibly even the empty face.  An important contribution of the
algebraic approach here is the identification of the minimal set of
group relaxations needed to solve {\em all} programs in the family
$IP_{A,c}$ and of the particular relaxations necessary to solve any
given program in the family.

\section{Arithmetic Degree}

For an associated set $\tau$ of $IP_{A,c}$ there are only finitely
many standard pairs of $\mathcal O_c$ indexed by $\tau$ since there
are only finitely many standard polytopes of the form $Q_u^{\bar
  \tau}$. Borrowing terminology from \cite{STV}, we call the number of
standard pairs of the form $(\cdot, \tau)$ the {\em multiplicity} of
$\tau$ in $\mathcal O_c$ (abbreviated as $mult(\tau)$). The total
number of standard pairs of $\mathcal O_c$ is called the {\em
  arithmetic degree} of $\mathcal O_c$. Our main goal in this section
is to provide bounds for these invariants of the the family $IP_{A,c}$
and discuss their relevance. We will need the following interpretation
from Section~3.

\begin{corollary}
  The multiplicity of the face $\tau$ of $\Delta_c$ in $\mathcal O_c$
  is the number of distinct standard polytopes of $IP_{A,c}$ indexed by 
  ${\bar \tau}$, and the arithmetic degree of
  $\mathcal O_c$ is the total number of standard polytopes of
  $IP_{A,c}$. 
\end{corollary}

\begin{proof}
This result follows from Theorem~\ref{stdpairs-equivs}.
\end{proof}

\noindent{\bf Example~\ref{knapsack} continued.}
The multiplicity of the associated set $\{3\}$ is eight while the
empty set has multiplicity twelve. The arithmetic degree of $\mathcal
O_c$ is hence twenty. \qed\\

If the standard pair decomposition of $\mathcal O_c$ is known, then we
can solve all programs in $IP_{A,c}$ by solving (arithmetic
degree)-many {\em linear systems} as follows. For a given $b \in
\mathbb N A$ and a standard pair $(u, \tau)$, consider the linear
system
\begin{equation}\label{linsys}
A_{\bar \tau}\pi_{\tau}(u) + A_{\tau}x = b, \,\,{\text or \,\,
  equivalently,} \,\, A_{\tau}x = b - A_{\bar \tau}\pi_{\tau}(u).
\end{equation}
As $\tau$ is a face of $\Delta_c$, this linear system can be solved
uniquely for $x$. Since the optimal solution of $IP_{A,c}(b)$ lies in
$S(w,\sigma)$ for some standard pair $(w, \sigma)$ of $\mathcal O_c$,
at least one non-negative and integral solution for $x$ will be found
as we solve the linear systems (\ref{linsys}) obtained by varying
$(u,\tau)$ over all the standard pairs of $\mathcal O_c$. If the
standard pair $(u,\tau)$ yields such a solution $v$, then
$(\pi_{\tau}(u),v)$ is the optimal solution of $IP_{A,c}(b)$. This
pre-processing of $IP_{A,c}$ has the same flavor as \cite{Kan93}. The
main result in \cite{Kan93} is that given a coefficient matrix $A \in
\mathbb R^{m \times n}$ and cost vector $c$, there exists {\em floor
  functions} $f_1, \ldots, f_k : \mathbb R^m \rightarrow \mathbb Z^n$
such that for a right hand side vector $b$, the optimal solution of
the corresponding integer program is the one among $f_1(b), \ldots,
f_k(b)$ that is feasible and attains the best objective function
value. The crucial point is that this algorithm runs in time bounded
above by a polynomial in the length of the data for fixed $n$ and $j$,
where $j$ is the affine dimension of the space of right hand sides.
Given this result, it is interesting to bound arithmetic degree.

The second equation in (\ref{linsys}) suggests that one could think of
the first arguments $u$ in the standard pairs $(u,\tau)$ of $\mathcal
O_c$ as ``correction vectors'' that need to be applied to find the
optimal solutions of programs in $IP_{A,c}$. Thus the arithmetic
degree of $\mathcal O_c$ is the total number of correction vectors
that are needed to solve all programs in $IP_{A,c}$. The
multiplicities of associated sets give a finer count of these
correction vectors, organized by faces of $\Delta_c$. If the optimal
solution of $IP_{A,c}(b)$ lies in the affine semigroup $S(w,\sigma)$
given by the standard pair $(w,\sigma)$ of $\mathcal O_c$, then $w$ is
a correction vector for this $b$ as well as all other $b$'s in $(Aw +
\mathbb N A_{\sigma})$. One obtains all correction vectors for
$IP_{A,c}$ by solving the (arithmetic degree)-many integer programs
with right hand sides $Au$ for all standard pairs $(u,\tau)$ of
$\mathcal O_c$. See \cite{Wol81} for a similar result from the
classical theory of group relaxations.

In Example~\ref{knapsack}, $\Delta_c = \{\{3\}\}$ and both its faces
$\{3\}$ and $\emptyset$ are associated to $IP_{A,c}$. In general, not
all faces of $\Delta_c$ need be associated sets of $IP_{A,c}$ and the
poset of associated sets can be quite complicated. (We will study this
poset in Section~5.) Hence, for $\tau \in \Delta_c$,
$mult(\tau) = 0$ unless $\tau$ is an associated set of $IP_{A,c}$. We
will now prove that all maximal faces of $\Delta_c$ are associated
sets of $IP_{A,c}$. Further, if $\sigma$ is a maximal face of
$\Delta_c$ then $mult(\sigma)$ is the absolute value of
$det(A_{\sigma})$ divided by the g.c.d. of the maximal minors of $A$.
This g.c.d is non-zero since $A$ has full row rank. If the columns of
$A$ span an affine hyperplane, then the absolute value of
$det(A_{\sigma})$ divided by the g.c.d. of the maximal minors of $A$
is called the {\em normalized volume} of the face $\sigma$ in
$\Delta_c$. We first give a non-trivial example.

\begin{example} \label{long-chain}
Consider the rank three matrix 
$$ A = \left[ \begin{array}{cccccc} 
5 & 0 & 0 & 2 & 1 & 0 \\
0 & 5 & 0 & 1 & 4 & 2 \\
0 & 0 & 5 & 2 & 0 & 3
\end{array} \right]$$
and the generic cost vector $c = (21, \, 6, \, 1, \, 0, \, 0, \, 0)$.
The first three columns of $A$ generate $cone(A)$ which is simplicial.
The regular triangulation $$\Delta_c = \{\{1, 3, 4\}, \, \{1,4,5\},
\,\{2,5,6\}, \, \{3,4,6\}, \,\{4,5,6\}\}$$
is shown in Figure~\ref{f6}
as a triangulation of $conv(A)$. The six columns of $A$ have been
labeled by their column indices.  The arithmetic degree of $\mathcal
O_c$ in this example is $70$. The following table shows all the
standard pairs organized by associated sets and the multiplicity of
each associated set. Note that all maximal faces of $\Delta_c$ are
associated to $IP_{A,c}$. The g.c.d.  of the maximal minors of $A$ is
five. Check that $mult(\sigma)$ is the normalized volume of $\sigma$
whenever $\sigma$ is a maximal face of $\Delta_c$.
\begin{figure}
\epsfig{file=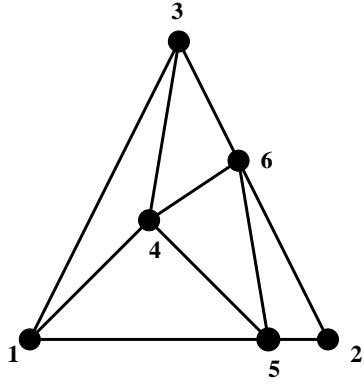, height=5cm}
\caption{The regular triangulation $\Delta_c$ for Example~\ref{long-chain}.}  
\label{f6}
\end{figure}
$$\begin{array}{|c|c|c|} \hline
  \tau & \text{ standard pairs } (\cdot, \tau) & mult(\tau) \\
  \hline \{1, 3, 4\} & (0, \cdot), \, (e_5, \cdot), \, (e_6, \cdot),
  \,
  (e_5+e_6, \cdot), \, (2e_6, \cdot) & 5\\
  \hline \{1,4,5\} & (0, \cdot),\, (e_2, \cdot),\, (e_3, \cdot),\,
  (e_6, \cdot), \,
  (e_2+e_3, \cdot),\, (2e_2, \cdot), & 8\\
  & (3e_2,\cdot),\, (2e_2+e_3, \cdot) &\\
  \hline
  \{2,5,6\}&  (0, \cdot), \, (e_3, \cdot), \, (2e_3, \cdot) & 3\\
  \hline \{3,4,6\}& (0, \cdot), \, (e_5, \cdot), \, (2e_5, \cdot), \,
  (3e_5, \cdot) & 4\\
  \hline \{4,5,6\}& (0, \cdot), \, (e_3, \cdot), \, (2e_3, \cdot), \,
  (3e_3, \cdot), \, (4e_3, \cdot) & 5\\
  \hline \{1,4\} & (e_3+2e_5+e_6, \cdot),\, (2e_3+2e_5+e_6, \cdot), \,
  (2e_3+2e_5, \cdot), & 5\\
  & (2e_3+3e_5, \cdot), \,  (2e_3+4e_5, \cdot) & \\
  \hline \{1,5\} & (e_2+e_6, \cdot), \, (2e_2+e_6, \cdot), \,
  (3e_2+e_6, \cdot)
  & 3\\
  \hline
  \{2,5\} & (e_3+e_4, \cdot), \, (e_4, \cdot), \, (2e_4, \cdot) & 3 \\
  \hline \{3,4\} & (e_2, \cdot), \, (e_1+e_2, \cdot), \, (e_1+2e_5,
  \cdot), \,
  (e_1+2e_5+e_6, \cdot), \, (e_2+e_5, \cdot), \,  & 5\\
  \hline
  \{3,6\} & (e_2, \cdot), \, (e_2 + e_5, \cdot) & 2\\
  \hline \{4,5\} & (e_2 + 2e_3, \cdot), \, (e_2+3e_3, \cdot), \,
  (2e_2+2e_3, \cdot), \, (3e_2+e_3, \cdot), \, (4e_2, \cdot) &
  5\\
  \hline
  \{5,6\} & (e_2+3e_3, \cdot) & 1\\
  \hline \{1\} & (e_2+e_3+e_6, \cdot), \, (e_2+e_3+e_5+e_6, \cdot), \,
  (e_2+2e_6, \cdot), & 6 \\
  & (e_2+e_3+2e_6, \cdot),
  (2e_2+2e_6, \cdot), \, (e_2+e_3+2e_5+e_6, \cdot) & \\
  \hline
  \{3\} & (e_1+e_2+e_6, \cdot), \, (e_1+e_2+2e_6, \cdot) & 2\\
  \hline \{4\} & (e_1+e_2+2e_3+e_5, \cdot), \,
  (e_1+e_2+2e_3+2e_5, \cdot), & 6 \\
  & (e_1+e_2+2e_3+3e_5, \cdot), \,
  (e_1+e_2+2e_3+4e_5, \cdot), & \\
  & (e_1+3e_3+3e_5, \cdot), \,
  (e_1+3e_3+4e_5, \cdot) & \\
  \hline \{\emptyset\} & (e_1+e_2+2e_3+e_5+e_6, \cdot), \,
  (e_1+e_2+2e_3+2e_5+e_6, \cdot), & 7 \\
  & (e_1+2e_2+e_3+e_6, \cdot), \,
  (e_1+2e_2+e_3+e_5+e_6, \cdot), & \\
  & (e_1+2e_2+e_3+2e_5+e_6, \cdot), \,
  (e_1+2e_2+e_3+2e_6, \cdot), & \\
  &  (e_1+3e_2+2e_6, \cdot) & \\
  \hline
  & \text{Arithmetic \,\,Degree} & 70 \\
  \hline
\end{array}$$\\

Observe that the integer program $IP_{A,c}(b)$ where $b =
A(e_1+e_2+e_3)$ is solved by $G^{\tau}(b)$ with $\tau = \{1, 4, 5\}$.
By Proposition~\ref{optbases}, Gomory's relaxation of $IP_{A,c}(b)$ is
indexed by $\sigma = \{4,5,6\}$ since $b$ lies in the interior of the
face $cone(A_{\sigma})$ of $\Delta_c$. However, neither this
relaxation nor any nontrivial extended relaxation solves $IP_{A,c}(b)$
since the optimal solution $e_1+e_2+e_3$ is not covered by any
standard pair $(\cdot, \tau)$ where $\tau$ is a non-empty subset of
$\{4,5,6\}$. 
\qed\\
\end{example}

\begin{theorem} \label{maxsimplex}
For a set $\sigma \subseteq \{1, \ldots, n\}$, $(0,\sigma)$ is a
standard pair of $\mathcal O_c$ if and only if $\sigma$ is a maximal
face of $\Delta_c$.
\end{theorem}

\begin{proof} 
  If $\sigma$ is a maximal face of $\Delta_c$, then by
  Definition~\ref{regtriang}, there exists $y \in \mathbb R^d$ such
  that $yA_{\sigma} = c_{\sigma}$ and $yA_{\bar \sigma} < c_{\bar
    \sigma}$.  Then $p = c_{\bar \sigma} - yA_{\bar \sigma} > 0$ and
  $pB^{\bar \sigma} = (c_{\bar \sigma} - yA_{\bar \sigma})B^{\bar
    \sigma} = c_{\bar \sigma}B^{\bar \sigma} - yA_{\bar \sigma}B^{\bar
    \sigma} = c_{\bar \sigma}B^{\bar \sigma} +yA_{\sigma}B^{\sigma} =
  c_{\bar \sigma}B^{\bar \sigma} + c_{\sigma}B^{\sigma} = cB.$ Hence
  there is a positive dependence relation among $(-cB)$ and the rows
  of $B^{\bar \sigma}$. Since $\sigma$ is a maximal face of
  $\Delta_c$, $|det(A_{\sigma})| \neq 0$. However, $|det(B^{\bar
    \sigma})| = |det(A_{\sigma})|$ which implies that $|det(B^{\bar
    \sigma})| \neq 0$. Therefore, $(-cB)$ and the rows of $B^{\bar
    \sigma}$ span $\mathbb R^{n-d}$ positively. This implies that
  $Q_0^{\bar \sigma} = \{ z \in \mathbb R^{n-d} \,\, : \,\, B^{\bar
    \sigma}z \leq 0, \,\, (-cB) \cdot z \leq 0 \}$ is a polytope
  consisting of just the origin. If any inequality defining this
  simplex is dropped, the resulting relaxation is unbounded as only
  $n-d$ inequalities would remain. Hence $Q_0^{\bar \sigma}$ is a
  standard polytope of $IP_{A,c}$ and by
  Theorem~\ref{stdpairs-equivs}, $(0,\sigma)$ is a standard pair of
  $\mathcal O_c$.

Conversely, if $(0,\sigma)$ is a standard pair of $\mathcal O_c$ then
$Q_0^{\bar \sigma}$ is a standard polytope of $IP_{A,c}$. Since every
inequality in the definition of $Q_0^{\bar \sigma}$ gives a halfspace
containing the origin and $Q_0^{\bar \sigma}$ is a polytope,
$Q_0^{\bar \sigma} = \{0\}$. Hence there is a positive linear
dependence relation among $(-cB)$ and the rows of $B^{\sigma}$. If
$|\bar \sigma| > n-d$, then $Q_0^{\bar \sigma}$ would coincide with
the relaxation obtained by dropping some inequality from those in
$B^{\bar \sigma}z \leq 0$. This would contradict that $Q_0^{\bar
  \sigma}$ was a standard polytope and hence $|\sigma| = d$ and
$\sigma$ is a maximal face of $\Delta_c$.
\end{proof}

\begin{corollary} Every maximal face of $\Delta_c$ is an associated
  set of $IP_{A,c}$. 
\end{corollary}

For Theorem~\ref{maxmultiplicity} and Corollary~\ref{lowerbound} below
we assume that the g.c.d. of the maximal minors of $A$ is one which
implies that $\mathbb Z A = \mathbb Z^d$.

\begin{theorem} \label{maxmultiplicity}
If $\sigma$ is a maximal face of $\Delta_c$ then 
the multiplicity of $\sigma$ in $\mathcal O_c$ is $|det(A_{\sigma})|$.
\end{theorem}

\begin{proof}
  Consider the full dimensional lattice $\mathcal L_{\sigma} =
  \pi_{\sigma}(\mathcal L) = \{ B^{\bar \sigma}z \,\,:\,\, z \in
  \mathbb Z^{n-d} \}$ in $\mathbb Z^{n-d}$. Since the g.c.d. of the
  maximal minors of $A$ is assumed to be one, the lattice $\mathcal
  L_{\sigma}$ has index $|det(B^{\bar \sigma})| = |det(A_{\sigma})|$
  in $\mathbb Z^{n-d}$. Since $\mathcal L_{\sigma}$ is full
  dimensional, it has a strictly positive element which guarantees
  that each equivalence class of $\mathbb Z^{n-d}$ modulo $\mathcal
  L_{\sigma}$ has a non-negative member. This implies that there are
  $|det(A_{\sigma})|$ distinct equivalence classes of $\mathbb
  N^{n-d}$ modulo $\mathcal L_{\sigma}$. Recall that if $u$ is a
  feasible solution to $IP_{A,c}(b)$ then 
  $$G^{\sigma}(b) = minimize \,\,\{{\tilde c_{\bar \sigma}} \cdot
  x_{\bar \sigma} \,\,:\,\, x_{\bar \sigma} \equiv u_{\bar \sigma}
  \,\,(mod \,\, \mathcal L_{\sigma}),\,\, x_{\bar \sigma} \in \mathbb
  N^{n-d} \}.$$
  Since there are $|det(A_{\sigma})|$ equivalence classes of $\mathbb
  N^{n-d}$ modulo $\mathcal L_{\sigma}$, there are $|det(A_{\sigma})|$
  distinct group relaxations indexed by $\sigma$. The optimal solution
  of each program becomes the right hand side vector of a standard
  polytope (simplex) of $IP_{A,c}$ indexed by $\sigma$. Since no two
  optimal solutions are the same (as they come from different
  equivalence classes of $\mathbb N^{n-d}$ modulo $\mathcal
  L_{\sigma}$), there are precisely $|det(A_{\sigma})|$ standard
  polytopes of $IP_{A,c}$ indexed by $\sigma$.
\end{proof}

\begin{corollary} \label{lowerbound}
  The arithmetic degree of $\mathcal O_c$ is bounded
  below by the sum of the absolute values of $det(A_{\sigma})$ as
  $\sigma$ varies among the maximal faces of $\Delta_c$.
\end{corollary}

A {\em primary} ideal $J$ in $k[x_1, \ldots, x_n]$ is a proper ideal
such that $fg \in J$ implies either $f \in J$ or $g^t \in J$ for some
positive integer $t$. A {\em prime} ideal $J$ of $k[x_1, \ldots, x_n]$
is a proper ideal such that $fg \in J$ implies that either $f \in J$
or $g \in J$. A {\em primary decomposition} of an ideal $I$ in $k[x_1,
\ldots, x_n]$ is an expression of $I$ as a finite intersection of {\em
  primary ideals} in $k[x_1, \ldots, x_n]$. Lemma~3.3 in \cite{STV}
shows that every monomial ideal $M$ in $k[x_1, \ldots, x_n]$ admits a
primary decomposition into irreducible primary ideals that are indexed
by the standard pairs of $M$. The {\em radical} of an ideal $I \subset
k[x_1, \ldots, x_n]$ is the ideal $\sqrt{I} := \{ f \in S \,\, : \,\,
f^t \in I, \,\, {\text for \,\,some\,\, positive\,\, integer\,\, t}
\}$. Radicals of primary ideals are prime. The radicals of the primary
ideals in a minimal primary decomposition of an ideal $I$ are called
the {\em associated primes} of $I$. This list of prime ideals is
independent of the primary decomposition of the ideal. The minimal
elements among the associated primes of $I$ are called the {\em
  minimal primes} of $I$ while the others are called the {\em embedded
  primes} of $I$. The minimal primes of $I$ are precisely the defining
ideals of the isolated components of the zero-set or {\em variety} of
$I$ while the embedded primes cut out embedded subvarieties in the
isolated components. See a textbook in commutative algebra
like \cite{Eis} for more details.

A face $\tau$ of $\Delta_c$ is an associated set of $IP_{A,c}$ if and
only if the monomial prime ideal $p_{\tau} := \langle x_j \,\, : \,\,
j \not \in \tau \rangle$ is an associated prime of the ideal
$in_c(I_A)$. Further, $p_{\sigma}$ is a minimal prime of $in_c(I_A)$
if and only if $\sigma$ is a maximal face of $\Delta_c$. Hence the
lower dimensional associated sets of $IP_{A,c}$ index the embedded
primes of $in_c(I_A)$. The standard pair decomposition of a monomial
ideal was introduced in \cite{STV} to study its associated primes. The
multiplicity of an associated prime $p_{\tau}$ of $in_c(I_A)$ is an
algebraic invariant of $in_c(I_A)$, and \cite{STV} shows that this is
exactly the number of standard pairs indexed by $\tau$. Similarly, the
arithmetic degree of $in_c(I_A)$ is a refinement of the geometric
notion of {\em degree} and \cite{STV} shows that this number is the
total number of standard pairs of $in_c(I_A)$. These connections
explain our choice of terminology.  Theorem~\ref{maxsimplex} is a
translation of the specialization of Lemma~3.5 in \cite{STV} to toric
initial ideals. We refer the interested reader to \cite[\S 8 and \S
12.D]{Stu} and \cite[\S 3]{STV} for the algebraic connections.
Theorem~\ref{maxmultiplicity} is a staple result of toric geometry and
also follows from \cite[Theorem 1]{Gom65}. It is proved via 
the algebraic technique of {\em localization} in \cite[Theorem 8.8]{Stu}. 

Theorem~\ref{maxmultiplicity} gives a precise bound on the
multiplicity of a maximal associated set of $IP_{A,c}$, which in turn
provides a lower bound for the arithmetic degree of $\mathcal O_c$ in
Corollary~\ref{lowerbound}. No exact result like
Theorem~\ref{maxmultiplicity} is known when $\tau$ is a lower
dimensional associated set of $IP_{A,c}$. Such bounds would provide a
bound for the arithmetic degree of $\mathcal O_c$. Bounds on the
arithmetic degree of a general monomial ideal in terms of its
dimension and minimal generators can be found in \cite[Theorem
3.1]{STV}. One hopes that stronger bounds are possible for toric
initial ideals. We close with a first attempt at bounding the
arithmetic degree of $\mathcal O_c$ (under certain non-degeneracy
assumptions). This result is due to Ravi Kannan, and its simple
arguments are along the lines of proofs in \cite{Kan92}
and \cite{KLS}. 

Suppose $S \in \mathbb Z^{m \times n}$ and $u \in \mathbb N^m$ are
fixed and $K_u := \{ x \in \mathbb R^n \,\, : \,\, Sx \leq u \}$ is
such that $K_u \cap \mathbb Z^n = \{0\}$ and the removal of any
inequality defining $K_u$ will bring in a non-zero lattice point into
the relaxation. Let $s^{(i)}$ denote the $i$th row of $S$, $M :=
max ||s^{(i)}||_1$ and $\Delta_k(S)$ and $\delta_k(S)$ be the
maximum and minimum absolute values of the $k \times k$
subdeterminants of $S$. We will assume that $\delta_n(S) \neq 0$ which
is a non-degeneracy condition on the data. We assume this set up
in Theorem~\ref{kannan} and Lemmas~\ref{a} and \ref{b}.

\begin{definition} If $K$ is a convex set and $v$ a non-zero vector in 
$\mathbb R^n$, the {\bf width of $K$ along $v$}, denoted as
${\mathit width}_v(K)$ is $max \,\,\{ v \cdot x \,\,:\,\, x \in K \} - 
min\,\,\{ v \cdot x \,\,:\,\, x \in K \}$.
\end{definition}

Note that ${\mathit width}_v(K)$ is invariant under translations of $K$.

\begin{theorem} \label{kannan}
If $K_u$ is as above then 
$ 0 \leq u_i \leq 2M(n+2) \frac{\Delta_n(S)}{\delta_n(S)}.$
\end{theorem}

\begin{lemma} \label{a}
  If $K_u$ is as above then for some $t$, $1 \leq t \leq m$, 
  $ {\mathit width}_{s^{(t)}}(K_u) \leq M(n+2).$
\end{lemma}

\begin{proof} Clearly, $K_u$ is bounded since otherwise there would be
  a non-zero lattice point on an unbounded edge of $K_u$ due to the
  integrality of all data. Suppose ${\mathit width}_{s^{(t)}}(K_u) > M(n+2)$ for
  all rows $s^{(t)}$ of $S$. Let $p$ be the center of gravity of
  $K_u$. Then by a property of the center of gravity, for any $x \in
  K_u$, $(1/(n+1))$th of the vector from $p$ to the reflection of $x$
  about $p$ is also in $K_u$, i.e., $(1 + \frac{1}{n+1})p -
  \frac{1}{n+1}x \in K_u$. Fix $i$, $1 \leq i \leq m$ and let $x_0$
  minimize $s^{(i)} \cdot x$ over $K_u$. By the definition of width,
  we then have $u_i - s^{(i)} \cdot x_0 > M(n+2)$ which implies that
\begin{equation} \label{one}
s^{(i)} \cdot x_0 < u_i - M(n+2).
\end{equation}
Now $s^{(i)} \cdot ((1 + \frac{1}{n+1})p - \frac{1}{n+1}x_0) \leq
u_i$ implies that 
\begin{equation} \label{two}
s^{(i)} \cdot p \leq u_i (\frac{n+1}{n+2}) + \frac{s^{(i)} \cdot x_0}{n+2}.
\end{equation}
Combining (\ref{one}) and (\ref{two}) we get 
\begin{equation} \label{three}
s^{(i)} \cdot p < u_i - M.
\end{equation}

Let $q = \lfloor p \rfloor$ be the vector obtained by rounding down
all components of $p$. Then $p = q + r$ where $0 \leq r_j < 1$ for all
$j = 1, \ldots, n$, and by (\ref{three}), $s^{(i)} \cdot (q
+ r) < u_i - M$ which leads to $s^{(i)} \cdot q + (s^{(i)} \cdot r +
M) < u_i$. Since $M = max ||s^{(i)}||_1$, 
\begin{equation} \label{four}
-M \leq s^{(i)} \cdot r \leq M.
\end{equation}
and hence, $s^{(i)} \cdot q < u_i$. Repeating this argument for all rows
of $S$, we get that $q \in K_u$. Similarly, if $q' = \lceil p \rceil$
is the vector obtained by rounding up all components of $p$, then $p =
q'-r$ where $0 \leq r_j < 1$ for all $j = 1, \ldots, n$. Then 
(\ref{three}) implies that $s^{(i)} \cdot (q'-r) < u_i - M$ which leads
to $s^{(i)} \cdot q' + (M - s^{(i)} \cdot r) < u_i$. Again by (\ref{four}), 
$s^{(i)} \cdot q' < u_i$ and hence $q' \in K_u$. 
Since $q \neq q'$, at least one of them is non-zero which contradicts
that $K_u \cap \mathbb Z^n = \{0\}$. 
\end{proof}

\begin{lemma} \label{b}
For any two rows $s^{(i)}, s^{(j)}$ of $S$, ${\mathit width}_{s^{(i)}}(K_u)
\leq 2 \frac{\Delta_n(S)}{\delta_n(S)} {\mathit width}_{s^{(j)}}(K_u)$.
\end{lemma}

\begin{proof}
  Without loss of generality we may assume that $j = n+1$. Since $K_u$
  is bounded, ${\mathit width}_{s^{(j)}}(K_u)$ is finite. Suppose the
  minimum of $s^{(j)} \cdot x$ over $K_u$ is attained at $v$. Since
  translations leave the quantities in the lemma invariant, we may
  prove the lemma for the body $K_{u'}$ obtained by translating $K_u$
  by $-v$. Now $s^{(j)} \cdot x$ is minimized over $K_{u'}$ at the
  origin. By LP duality, there are $n$ linearly independent
  constraints among the $m$ defining $K_{u'}$ such that the minimum of
  $s^{(n+1)}\cdot x$ subject to just these $n$ constraints is attained
  at 0. After renumbering the inequalities if necessary, assume these
  $n$ constraints are the first $n$. Let
  $$D = \{ x \,\,:\,\, s^{(l)} \cdot x \leq u'_l,\,\, l = 1,2,\ldots,
    n+1 \}$$
  where of course $u_1'=u_2'=\cdots=u_n'=0$. Then by the above, $D$ is
  a bounded simplex. 
  
  Since $D$ contains $K_{u'}$, it suffices to show that for each $i$,
\begin{equation} \label{five}
{\mathit width}_{s^{(i)}}(D) \leq 2 (\frac{\Delta_n(S)}{\delta_n(S)})
{\mathit width}_{s^{(n+1)}}(K_{u'}) = 2 (\frac{\Delta_n(S)}{\delta_n(S)})
u'_{n+1}.
\end{equation}
We show that for each vertex $q$ of $D$, $|s^{(i)} \cdot q | \leq
(\frac{\Delta_n(S)}{\delta_n(S)}) u'_{n+1}$ which will prove
(\ref{five}). This is clearly true for $q = 0$. Without loss of
generality assume that vertex $q$ satisfies 
$s^{(l)} \cdot q = u'_l$ for $l = 2,3,\ldots, n+1$. Since the
determinant of the submatrix of $S$ consisting of the rows $s^{(2)},
\ldots, s^{(n+1)}$ is not zero, for any $i$ there
exists rationals $\lambda_l$ such that $s^{(i)} = \sum_{l=2}^{n+1}
\lambda_l s^{(l)}$. By Cramer's rule, $|\lambda_l| \leq
(\frac{\Delta_n(S)}{\delta_n(S)})$. Therefore, $s^{(i)} \cdot q  = 
\sum_{l=2}^{n+1} \lambda_l s^{(l)} \cdot q = \sum_{l=2}^{n+1}
\lambda_l u'_l = \lambda_{n+1} u'_{n+1}$ since $u'_l = 0$ for $l = 2,
\ldots, n$. This proves that 
$$ |s^{(i)} \cdot q| = |\lambda_{n+1} u'_{n+1}| =
|\lambda_{n+1}|u'_{n+1} \leq (\frac{\Delta_n(S)}{\delta_n(S)})
u'_{n+1}.$$
\end{proof}

\noindent{\em Proof of Theorem~\ref{kannan}.} 
From Lemmas~\ref{a} and \ref{b} it follows that for any $i$, $1 \leq i
\leq m$, ${\mathit width}_{s^{(i)}}(K_u) \leq 2
(\frac{\Delta_n(S)}{\delta_n(S)}) M (n+2) =
2M(n+2)(\frac{\Delta_n(S)}{\delta_n(S)})$.  Since $0 \in K_u$, $min
\{s^{(i)} \cdot x \,\,:\,\, x \in K_u \} \leq 0$ while $max \{s^{(i)}
\cdot x \,\,:\,\, x \in K_u \} = u_i$.  Therefore, $u_i = u_i - 0 \leq
{\mathit width}_{s^{(i)}}(K_u)$ and hence, $0 \leq u_i \leq
2M(n+2)(\frac{\Delta_n(S)}{\delta_n(S)})$ for all $1 \leq i \leq m$.
\qed\\

Reverting back to our set up, let $\bar B = \left[ \begin{array}{cc} B
    \\ -cB \end{array} \right]$. Suppose $K_u$ is the standard
polytope $Q_u^{\bar \tau}$.  By Theorem~\ref{kannan}, $0 \leq u_i \leq
2M(n-d+2)(\frac{\Delta_n(\bar B)}{\delta_n(\bar B)})$.

\begin{corollary} \label{arithdegbound} 
  If no maximal minor of $\bar B$ is zero, then the arithmetic degree
  of $\mathcal O_c$ is at most $ \left( 2M(n-d+2)(\frac{\Delta_n(\bar
      B)}{\delta_n(\bar B)}) \right)^n$.
\end{corollary}

The above arguments do not use the condition that the
removal of an inequality from $K_u$ will bring in a lattice point into
the relaxation. Further, the bound is independent of the number of
facets of $K_u$, and Corollary~\ref{arithdegbound} is straightforward.
Thus, further improvements may be possible with more effort. However,
apart from providing a bound for arithmetic degree, these proofs have
the nice feature that they build a bridge to techniques from the
geometry of numbers that have played a central role in theoretical
integer programming as seen in the work of Kannan, Lenstra, Lov{\'a}sz,
Scarf and others. See \cite{Lov} for a survey.

\section{The Chain Theorem}

We now examine the structure of the poset of associated sets of
$IP_{A,c}$ which we denote as $Assets(IP_{A,c})$. All elements of
$Assets(IP_{A,c})$ are faces of the regular triangulation $\Delta_c$
and the partial order is set
inclusion. Theorem~\ref{maxsimplex} provides a first result.  

\begin{corollary} The maximal elements of $Assets(IP_{A,c})$ are
 the maximal faces of $\Delta_c$.
\end{corollary}

\begin{figure}
\epsfig{file=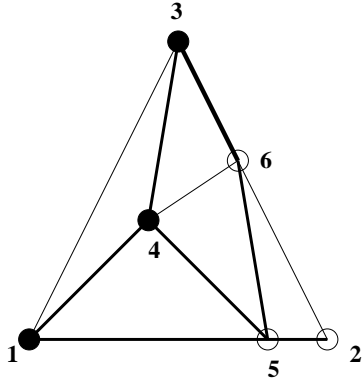, height=5cm}
\caption{ Lower dimensional associated sets of
  Example~\ref{long-chain} except the empty set.}
\label{f7}
\end{figure}

\noindent{\bf Example~\ref{long-chain} continued.} 
The lower dimensional associated sets of this example (except the
empty set) are the thick faces of $\Delta_c$ shown in Figure~\ref{f7}.
\qed\\

Despite the seemingly chaotic structure of $Assets(IP_{A,c})$ beyond
its maximal elements, it has an important structural property that we
now explain.

\begin{theorem} \label{chain-theorem} {\bf [The Chain Theorem]}
If $\tau \in \Delta_c$ is an associated set of $IP_{A,c}$ and $|\tau|
< d$ then there exists a face $\tau' \in \Delta_c$ that is also an
associated set of $IP_{A,c}$ with the property that $\tau \subset
\tau'$ and $|\tau' \backslash \tau| = 1$.
\end{theorem}

\begin{proof} Since $\tau$ is an associated set of $IP_{A,c}$,
  by Theorem~\ref{stdpairs-equivs}, $\mathcal O_c$ has a standard pair
  of the form $(v,\tau)$ and $Q_v^{\bar \tau} = \{z \in \mathbb
  R^{n-d} \, : \, B^{\bar \tau}z \leq \pi_{\tau}(v), \,\, (-cB) \cdot
  z \leq 0 \}$ is a standard polytope of $IP_{A,c}$. Since $|\tau| <
  d$, $\tau$ is not a maximal face of $\Delta_c$ and hence by
  Theorem~\ref{maxsimplex}, $v \neq 0$. For each $i \in \bar \tau$,
  let $R^i$ be the relaxation of $Q_v^{\bar \tau}$ obtained by
  removing the $i$th inequality $b_i \cdot z \leq v_i$ from $B^{\bar
    \tau}z \leq \pi_{\tau}(v)$, i.e., $$R^i := \{ z \in \mathbb
  R^{n-d} \, : B^{\bar \tau \backslash \{i\}}z \leq \pi_{\tau \cup
    \{i\}}(v), \, (-cB) \cdot z \leq 0\}.$$
  Let $E^i := R^i \backslash
  Q_v^{\bar \tau}$. Clearly, $E^i \cap Q_v^{\bar \tau} = \emptyset$,
  and, since the removal of $b_i \cdot z \leq v_i$ introduces at least
  one lattice point into $R^i$, $E^i \cap \mathbb Z^{n-d} \neq
  \emptyset$.  Let $z_i^*$ be the optimal solution to $
  minimize \,\, \{ (-cB) \cdot z : z \in E^i \cap 
  \mathbb Z^{n-d} \}$ if the program is bounded.  This integer program
  is always feasible since $E^i \cap \mathbb Z^{n-d} \neq \emptyset$,
  but it may not have a finite optimal value.  However, there exists
  at least one $i \in \bar \tau$ for which the above integer program
  is bounded.  To see this, pick a maximal simplex $\sigma \in
  \Delta_c$ such that $\tau \subset \sigma$. The polytope $\{z \in
  \mathbb R^{n-d}: B^{\bar \sigma}z \leq \pi_{\sigma}(v),\,\, (-cB)
  \cdot z \leq 0 \}$ is a simplex and hence bounded. This polytope
  contains all $E^i$ for $i \in \sigma \backslash \tau$, and hence all
  these $E^i$ are bounded and have finite optima with respect to
  $(-cB) \cdot z$.  We may assume that the inequalities in $B^{\bar \tau}z
  \leq \pi_{\tau}(v)$ are labeled so that the finite optimal values
  are ordered as $(-cB) \cdot z _1^* \geq (-cB) \cdot z _2^* \geq
  \cdots \geq (-cB) \cdot z_p^*$ where $ \{1,2,\ldots,p\} \subseteq
  \bar \tau$.

\vskip 0.2cm 
\noindent{\bf Claim:} {\em Let $N^1 := \{z \in \mathbb R^{n-d} \, : \,
  B^{\bar \tau \backslash \{1\}}z \leq \pi_{{\tau \cup \{1\}}}(v),
  \,\,(-cB) \cdot z \leq (-cB) \cdot z_1^*\}$. Then $z_1^*$ is the
  unique lattice point in $N^1$ and the removal of any inequality from
  $B^{\bar \tau \backslash \{1\}}z \leq \pi_{{\tau \cup \{1\}}}(v)$
  will bring in a new lattice point into the relaxation.}
 
\noindent{\em Proof.} Since $z_1^*$ lies in $R^1$, $0 = (-cB) \cdot 0 \geq
(-cB) \cdot z_1^*$. However, $0 > (-cB) \cdot z_1^*$ since otherwise,
both $z_1^*$ and $0$ would be optimal solutions to $minimize \{(-cB)
\cdot z : z \in R^1 \}$ contradicting that $c$ is
generic. Therefore,\\
$\begin{array}{rl} N^1 = & R^1 \cap \{z \in \mathbb R^{n-d} \,\,:
  \,\,(-cB) \cdot z
  \leq (-cB) \cdot z _1^*\} \\
  = & (E^1 \cup Q_v^{\bar \tau}) \cap \{z \in \mathbb R^{n-d} \,\,:
  \,\, (-cB) \cdot z \leq (-cB) \cdot z_1^*\} \\
  = & (E^1 \cap \{z \in \mathbb R^{n-d} \,\,:\,\, (-cB) \cdot z \leq
  (-cB) \cdot z_1^*\}) \,\,\\ 
     & \quad \bigcup \,\, (Q_v^{\bar \tau} \cap \{z \in \mathbb R^{n-d} :
  (-cB) \cdot z \leq (-cB) \cdot z_1^*\}).
  \end{array}$\\
  Since $c$ is generic, $z_1^*$ is the unique lattice point in the
  first polytope and the second polytope is free of lattice points.
  Hence $z_1^*$ is the unique lattice point in $N^1$. The relaxation
  of $N^1$ got by removing $b_j \cdot z \leq v_j$ is the polyhedron
  $N^1 \cup (E^j \cap \{z \in \mathbb R^{n-d} \,\, : \,\, (-cB) \cdot
  z \leq (-cB) \cdot z_1^* \})$ for $j \in \bar \tau$ and $j \neq 1$.
  Either this is unbounded, in which case there is a lattice point $z$
  in this relaxation such that $(-cB) \cdot z_1^* \geq (-cB) \cdot z$,
  or (if $j \leq p$) we have $(-cB) \cdot z_1^* \geq (-cB) \cdot
  z_j^*$ and $z_j^*$ lies in this relaxation. \hfill $\diamondsuit$\\
  
  Translating $N^1$ by $-z_1^*$ we get $Q_{v'}^{\bar \tau \backslash
    \{1\}} := \{z \in \mathbb R^{n-d} : (-cB) \cdot z \leq 0,\,
  B^{\bar \tau \backslash \{1\}}z \leq v' \}$ where $v' = \pi_{\tau
    \cup \{1\}}(v) - B^{\bar \tau \backslash \{1\}}z_1^* \geq 0$ since
  $z_1^*$ is feasible for all inequalities except the first one.  Now
  $Q_{v'}^{\bar \tau \backslash \{1\}} \cap \mathbb Z^{n-d} = \{0\}$,
  and hence $(v',\tau \cup \{1\})$ is a standard pair of $\mathcal
  O_c$.
\end{proof}

\noindent{\bf Example~\ref{long-chain} continued.} 
The empty set is associated to $IP_{A,c}$ and $\emptyset \subset 
\{1\} \subset \{1,4\} \subset \{1,4,5\}$ is a saturated chain in
$Assets(IP_{A,c})$ that starts at the empty set. \qed \\

In algebraic language, the chain theorem says that the associated
primes of $in_c(I_A)$ occur in saturated chains. This was proved in
\cite[Theorem~3.1]{HT2}. When the cost vector $c$ is not generic,
$in_c(I_A)$ is no longer a monomial ideal, and its associated primes
need not come in saturated chains. See \cite[Remark 3.3]{HT2} for such
an example. An important open question in the algebraic study of
integer programming is to characterize all monomial ideals that can
appear as the initial ideal (with respect to some generic cost vector)
of a toric ideal. In our set up this amounts to characterizing all
down sets in $\mathbb N^n$ that can appear as the set of optimal
solutions to a family $IP_{A,c}$, where $A$ and $c$ satisfy the
assumptions from Section~2. Theorem~\ref{chain-theorem} imposes the
necessary condition that the poset of sets indexing the standard pairs
of the down set have the chain property. Unfortunately, this is not
sufficient to characterize down sets of the form $\mathcal O_c$. See
\cite{MSY} for another class of monomial ideals that also have the
chain property.

Since the elements of $Assets(IP_{A,c})$ are faces of $\Delta_c$, a
maximal face of which is a $d$-element set, the length of a maximal
chain in $Assets(IP_{A,c})$ is at most $d$. We denote the maximal
length of a chain in $Assets(IP_{A,c})$ by
$length(Assets(IP_{A,c}))$. When $n-d$ (the corank of $A$) is small
compared to $d$, $length(Assets(IP_{A,c}))$ has a stronger upper bound
than $d$. We use the following result of Bell and Scarf to prove the
bound.

\begin{theorem} \label{Bell-Scarf} \cite[Corollary 16.5 a]{Sch}
Let $Ax \leq b$ be a system of linear inequalities in $n$ variables, and 
let $c \in \mathbb R^n$. If max $\{ c \cdot x : Ax \leq b, \, x \in 
\mathbb Z^n \}$ is a finite number, then max $\{ c \cdot x : Ax \leq
b, \, x \in \mathbb Z^n \}$ =  max $\{  c \cdot x : A'x \leq b', x \in 
\mathbb Z^n \}$ for some subsystem $A'x \leq b'$ of $Ax \leq b$ with
at most $2^n - 1 $ inequalities. 
\end{theorem}

\begin{theorem} \label{length}
The length of a maximal chain in the poset of associated sets of 
$IP_{A,c}$ is at most ${\text min} (d,  2^{n-d} - (n-d+1))$. 
\end{theorem}

\begin{proof}
  As seen earlier, $length(Assets(IP_{A,c})) \leq d$. If $v$ lies in
  $\mathcal O_c$, then the origin is the optimal solution to the
  integer program $minimize \,\, \{(-cB) \cdot z \,\, : \,\, Bz \leq
  v, \,\,z \in \mathbb Z^{n-d} \}.$ By Theorem \ref{Bell-Scarf}, we
  need at most $2^{n-d}-1$ inequalities to describe the same integer
  program which means that we can remove at least $n - (2^{n-d} - 1)$
  inequalities from $Bz \leq v$ without changing the optimum.
  Assuming that the inequalities removed are indexed by $\tau$,
  $Q_v^{\bar \tau}$ will be a standard polytope of $IP_{A,c}$.
  Therefore, $|\tau| \geq n - (2^{n-d} - 1)$. This implies that the
  maximal length of a chain in $Assets(IP_{A,c})$ is at most $d - (n -
  (2^{n-d} - 1)) = 2^{n-d} - (n-d+1)$.
\end{proof}

\begin{corollary}
The cardinality of an associated set of $IP_{A,c}$ is at least 
$max(0,n - (2^{n-d} - 1))$.
\end{corollary}

\begin{corollary} \label{corank2-chain}
If $n-d = 2$, then $length(Assets(IP_{A,c})) \leq 1$.
\end{corollary}

\begin{proof} In this situation, $2^{n-d} - (n-d+1) = 4 - (4-2+1) =
4-3 = 1$.
\end{proof}

We conclude this section with a family of examples for which 
$length(Assets(IP_{A,c})) = 2^{n-d} - (n-d+1)$. This is adapted from 
\cite[Proposition~3.9]{HT2} which was modeled on a family of examples 
from \cite{PS1}.

\begin{proposition} \label{sharp}
  For each $m > 1$, there is an integer matrix $A$ of corank $m$ and a
  cost vector $c \in \mathbb Z^n$ where $n = 2^m-1$ such that
  $length(Assets(IP_{A,c})) = 2^m - (m+1)$.
\end{proposition} 

\begin{proof}
  Given $m > 1$, let $B' = (b_{ij}) \in \mathbb Z^{(2^m-1) \times m}$
  be the matrix whose rows are all the $\{1,-1\}$-vectors in ${\bf
    R}^m$ except $v = (-1,-1, \ldots, -1)$. Let $B \in \mathbb
  Z^{(2^m+m-1) \times m}$ be obtained by stacking $B'$ on top of
  $-I_m$ where $I_m$ is the $m \times m$ identity matrix. Set $n =
  2^m+m-1$, $d = 2^m-1$ and $A' = [I_d | B'] \in \mathbb Z^{d \times
    n}$. By construction, the columns of $B$ span the lattice $\{u \in
  \mathbb Z^n \,\,:\,\, A'u = 0 \}$. We may assume that the first row
  of $B'$ is $(1,1,\ldots,1) \in \mathbb R^m$. Adding this row to all
  other rows of $A'$ we get $A \in \mathbb N^{d \times n}$ with the
  same row space as $A'$. Hence the columns of $B$ are also a basis
  for the lattice $\{u \in \mathbb Z^n \,\,:\,\, Au = 0\}$. Since the
  rows of $B$ span $\mathbb Z^m$ as a lattice, we can find a cost
  vector $c \in \mathbb Z^n$ such that $(-cB) = v$. 
 
  For each row $b_i$ of $B'$ set $r_i := |\{b_{ij}\,\,:\,\,b_{ij} = 1
  \}|$, and let $r$ be the vector of all $r_i$s. By construction, the
  polytope $Q := \{z \in \mathbb R^m : B'z \leq r, -(cB) \cdot z \leq
  0 \}$ has no lattice points in its interior, and each of its $2^m$
  facets has exactly one vertex of the unit cube in $\mathbb R^m$ in
  its relative interior. If we let $w_i = r_i-1$, then the polytope
  $\{ z \in \mathbb R^m : B'z \leq w, -(cB) \cdot z \leq 0\}$ is a
  standard polytope $Q_{u}^{\bar \tau}$ of $IP_{A,c}$ where $\tau =
  \{d+1,d+2, \ldots, d+m = n\}$ and $w = \pi_{\tau}(u)$. Since a maximal
  face of $\Delta_c$ is a $d=(2^m - 1)$-element set and $|\tau| = m$,
  Theorem \ref{chain-theorem} implies that $length(Assets(IP_{A,c}))
  \geq 2^m - 1 - m = 2^m - (m+1)$. However, by Theorem~\ref{length}, 
   $length(Assets(IP_{A,c})) = min(2^m-1, 2^m - (m+1)) = 2^m - (m+1)$
  since $m > 1$ by assumption.
\end{proof}

\begin{example} If we choose $m = 3$ then $n = 2^m + m - 1 = 10$ and $d =
  2^m - 1 = 7$. Constructing $B'$ and $A$ as in
  Proposition~\ref{sharp}, we get 
$$ B' = \left[ \begin{array}{rrr} 
1 & 1 & 1 \\
-1 & 1 & 1 \\
1 & -1 & 1 \\
1 & 1 & -1 \\
-1 & -1 & 1 \\
-1 & 1 & -1 \\
1 & -1 & -1 \\
\end{array} \right] \,\,\,\,\, {\text and} \,\,\,\,\,
A = 
\left[ \begin{array}{cccccccccc} 
1 & 0 & 0 & 0 & 0 & 0 & 0 & 1 & 1 & 1 \\
1 & 1 & 0 & 0 & 0 & 0 & 0 & 0 & 2 & 2 \\
1 & 0 & 1 & 0 & 0 & 0 & 0 & 2 & 0 & 2 \\
1 & 0 & 0 & 1 & 0 & 0 & 0 & 2 & 2 & 0 \\
1 & 0 & 0 & 0 & 1 & 0 & 0 & 0 & 0 & 2 \\
1 & 0 & 0 & 0 & 0 & 1 & 0 & 0 & 2 & 0 \\
1 & 0 & 0 & 0 & 0 & 0 & 1 & 2 & 0 & 0 \\
\end{array} \right].$$
The vector $c = (11,0,0,0,0,0,0,10,10,10)$ satisfies $(-cB) =
(-1,-1,-1)$. The associated sets of $IP_{A,c}$ along with their 
multiplicities are given below.
\begin{center}
$\begin{array}{|c|c|c|c|}
\hline
{\tau} & {\text Multiplicity} & 
{\tau} & {\text Multiplicity}\\
\hline
\{4,5,6,7,8,9,10\}^* &    4& \{2,3,7,8,9,10\} &    2\\
\{1,5,6,7,8,9,10\} &    4& \{5,6,7,8,9,10\}^* &    1\\
\{3,4,6,7,8,9,10\} &    4& \{4,5,6,7,8,9\} &    1\\
\{2,3,4,6,7,9,10\} &    2& \{2,4,7,8,9,10\} &    2\\
\{2,3,4,7,8,9,10\} &    4& \{1,5,7,8,9,10\} &    1\\
\{3,4,5,6,7,8,10\} &    2& \{2,3,4,8,9,10\} &    1\\
\{2,3,4,5,6,7,10\} &    1& \{4,5,7,8,9,10\} &    2\\
\{2,4,5,6,7,9,10\} &    2& \{2,5,6,7,9,10\} &    1\\
\{2,3,6,7,9,10\} &    1&\{4,5,6,8,9,10\} &    2\\
\{3,4,5,6,8,10\} &    1&\{1,5,6,8,9,10\} &    1\\
\{2,4,5,7,9,10\} &    1&\{3,4,6,8,9,10\} &       2\\
\{1,6,7,8,9,10\} &    1&\{6,7,8,9,10\}^* &    1\\
\{3,5,6,7,8,10\} &    1&\{7,8,9,10\}^* & 1\\
\{3,6,7,8,9,10\} &    2&\{8,9,10\}^* &    1\\
\hline
\end{array}$
\end{center}
The elements in the unique maximal chain in $Assets(IP_{A,c})$ are 
marked with a $*$ and $length(Assets(IP_{A,c})) = 2^3 - (3+1) = 4$ as
predicted by Proposition~\ref{sharp}. 
\qed
\end{example}

\section{Gomory Integer Programs}

Recall from Definition~\ref{gomory-relaxations} that a group
relaxation $G^{\sigma}(b)$ of $IP_{A,c}(b)$ is called a Gomory
relaxation if $\sigma$ is a maximal face of $\Delta_c$. As discussed
in Section~2, these relaxations are the easiest to solve among all
relaxations of $IP_{A,c}(b)$. Hence it is natural to ask under what
conditions on $A$ and $c$ would all programs in $IP_{A,c}$ be solvable
by Gomory relaxations. We study this question in this section. The
majority of the results here are taken from \cite{HT01}.

\begin{definition} \label{gomory-family}
The family of integer programs $IP_{A,c}$ is a 
{\bf Gomory family} if, for {\em every} $b \in \mathbb N A$,
$IP_{A,c}(b)$ is solved by a group relaxation $G^{\sigma}(b)$
where $\sigma$ is a {\em maximal} face of the regular triangulation
$\Delta_c$. 
\end{definition}

\begin{theorem} \label{gomoryfam-theorem}
The following conditions are equivalent:\\
(i) $IP_{A,c}$ is a Gomory family. \\
(ii) The associated sets of $IP_{A,c}$ are precisely the maximal 
faces of $\Delta_c$.\\
(iii) $(\ast, \tau)$ is a standard pair of $\mathcal O_c$ if and only
if $\tau$ is a maximal face of $\Delta_c$.\\ 
(iv) All standard polytopes of $IP_{A,c}$ are simplices.
\end{theorem}

\begin{proof}
  By Definition~\ref{gomory-family}, $IP_{A,c}$ is a Gomory family if
  and only if for all $b \in \mathbb N A$, $IP_{A,c}(b)$ can be solved
  by one of its Gomory relaxations.  By Theorem~\ref{cover=solve},
  this is equivalent to saying that every $u \in \mathcal O_c$ lies in
  some $S(\ast, \sigma)$ where $\sigma$ is a maximal face of
  $\Delta_c$ and $(\ast, \sigma)$ a standard pair of $\mathcal O_c$.
  Definition~\ref{associatedset} then implies that all associated sets
  of $IP_{A,c}$ are maximal faces of $\Delta_c$. By
  Theorem~\ref{maxsimplex}, every maximal face of $\Delta_c$ is an
  associated set of $IP_{A,c}$ and hence $(i) \Leftrightarrow (ii)$.
  The equivalence of statements (ii), (iii) and (iv) follow from
  Theorem~\ref{stdpairs-equivs}.
\end{proof}

If $c$ is a generic cost vector such that for a triangulation $\Delta$
of $cone(A)$, $\Delta = \Delta_c$, then we say that $\Delta$ {\em
  supports} the order ideal $\mathcal O_c$ and the family of integer
programs $IP_{A,c}$. No regular triangulation of the matrix $A$ in
Example~\ref{long-chain} supports a Gomory family. Here is a matrix
with a Gomory family.

\begin{example} \label{example-with-gfamily}
Consider the $3 \times 6$ matrix 
$$A = \left[ \begin{array}{cccccc} 
1 & 0 & 1 & 1 & 1 & 1 \\
0 & 1 & 1 & 1 & 2 & 2 \\
0 & 0 & 1 & 2 & 3 & 4
\end{array} \right].$$
In this case, $cone(A)$ has 14 distinct regular triangulations and 48
distinct sets $\mathcal O_c$ as $c$ varies among all generic cost
vectors. Ten of these triangulations support Gomory families; one for
each triangulation. For
instance, if $c = 
(0,\,0,\,1,\,1,\,0,\,3)$, then
$$\Delta_c = 
\{ \sigma_1 = \{1, 2, 5\}, \,\, \sigma_2 = \{1, 4, 5\}, \,\, 
\sigma_3 = \{2,5,6\}, \,\, \sigma_4 = \{4,5,6\} \}$$
and $IP_{A,c}$ is a Gomory family since the standard pairs of
$\mathcal O_c$ are: $$(0, \sigma_1), \, (e_3, \sigma_1), \,  (e_4,
\sigma_1), \,  (0, \sigma_2), \, (0, \sigma_3), \, \text{and} \, (0,
\sigma_4).$$
\qed
\end{example}

Algebraically, $IP_{A,c}$ is a Gomory family if and only if the
initial ideal $in_c(I_A)$ has no embedded primes and hence
Theorem~\ref{gomoryfam-theorem} is a characterization of toric initial
ideals without embedded primes. A sufficient condition for an ideal in
$k[x_1, \ldots, x_n]$ to not have embedded primes is that it is {\em
  Cohen-Macaulay} \cite{Eis}. In general, Cohen-Macaulayness is not
necessary for an ideal to be free of embedded primes. However,
empirical evidence seemed to suggest for a while that for toric
initial ideals, Cohen-Macaulayness might be equivalent to being free
of embedded primes. A counterexample to this was found recently by
Laura Matusevich. The algebraic approach to integer programming allows
one to compute all down sets $\mathcal O_c$ of a fixed matrix $A$ as
$c$ varies among the set of generic cost vectors. See \cite{HuTh},
\cite{Stu} and \cite{ST} for details. The software package TiGERS
\cite{TiG} is custom-tailored for this purpose.

We now compare the notion of a Gomory family to the classical notion
of {\em total dual integrality} \cite[\S 22]{Sch}. It will be
convenient to assume that $\mathbb Z A = \mathbb Z^d$ for these
results.

\begin{definition} \label{tdi}
The system $yA \leq c$ is {\bf totally dual
integral (TDI)} if $LP_{A,c}(b)$ has an integral optimal  
solution for each $b \in cone(A) \cap \mathbb Z^d$. 
\end{definition}

\begin{definition}
The regular triangulation $\Delta_c$ is {\em unimodular} 
if $\mathbb{Z}A_\sigma = \mathbb Z^d$  for every maximal face 
$\sigma \in \Delta_c$. 
\end{definition}

\begin{example} The regular triangulation in
  Example~\ref{exs_regtriangs} (i) is unimodular while those in 
  Example~\ref{exs_regtriangs} (ii) and (iii) are not.  \qed\\
\end{example}

\begin{lemma} \label{TDI-ness} 
The system $yA \leq c$ is TDI if and only if 
the regular triangulation $\Delta_c$ is unimodular.
\end{lemma} 

\begin{proof} 
The regular triangulation $\Delta_c$ is the normal fan of $P_c$ by
Proposition~\ref{normalfan}, and it is unimodular if and only if 
$\mathbb Z A_{\sigma} = \mathbb Z^d$ for every
maximal face $\sigma \in \Delta_c$. This is equivalent to 
every $b \in cone(A_{\sigma}) \cap \mathbb Z^d$ lying in $\mathbb N
A_{\sigma}$ for every maximal face $\sigma$ of $\Delta_c$. By
Lemma~\ref{optbases}, this happens if and only if $LP_{A,c}(b)$
has an integral optimum for all $b \in cone(A) \cap \mathbb
Z^d$. 
\end{proof}

Corollary 8.4 in \cite{Stu} shows that $\Delta_c$ is unimodular if and
only if the monomial ideal $in_c(I_A)$ is generated by {\em
  square-free} monomials. Hence, by computing $in_c(I_A)$, one can
determine whether $yA \leq c$ is TDI. Such computations can
be carried out on computer algebra systems like CoCoA \cite{coc} or 
MACAULAY 2 \cite{GrS} for moderately sized examples. See \cite{Stu} for
algorithms. Standard pair decompositions of monomial ideals can
be computed with MACAULAY 2 \cite{HoSm}.

\begin{theorem}\label{TDI-Gomory}
If $yA \leq c$ is TDI then $IP_{A,c}$ is a Gomory
family.
\end{theorem}

\begin{proof}
  By Theorem~\ref{maxsimplex}, $(0, \sigma)$ is a standard pair of
  $\mathcal O_c$ for every maximal face $\sigma$ of $\Delta_c$.
  Lemma~\ref{TDI-ness} implies that $cone(A_{\sigma})$ is unimodular
  (i.e., $\mathbb Z A_{\sigma} = \mathbb Z^d$), and therefore $\mathbb
  N A_{\sigma} = cone(A_{\sigma}) \cap \mathbb Z^d$ for every maximal
  face $\sigma$ of $\Delta_c$.  Hence the semigroups $\mathbb N
  A_{\sigma}$ arising from the standard pairs $(0,\sigma)$ as $\sigma$
  varies over the maximal faces of $\Delta_c$ cover $\mathbb N A$.
  Therefore the only standard pairs of $\mathcal O_c$ are $(0,\sigma)$
  as $\sigma$ varies over the maximal faces of $\Delta_c$. The result
  then follows from Theorem~\ref{gomoryfam-theorem}.
\end{proof}

When $yA \leq c$ is TDI, the multiplicity of a maximal face
of $\Delta_c$ in $\mathcal O_c$ is one (from
Theorem~\ref{maxmultiplicity}). By Theorem~\ref{TDI-Gomory}, no lower
dimensional face of $\Delta_c$ is associated to $IP_{A,c}$.  While
this is sufficient for $IP_{A,c}(b)$ to be a Gomory family, it is far
from necessary. TDI-ness guarantees {\em local integrality} in the
sense that $LP_{A,c}(b)$ has an integral
optimum for every integral $b$ in $cone(A)$. In contrast, if
$IP_{A,c}$ is a Gomory family, the linear optima of the programs in
$LP_{A,c}$ may not be integral.

If $A$ is {\em unimodular} (i.e., $\mathbb Z A_{\sigma} = \mathbb Z^d$
for every nonsingular maximal submatrix $A_{\sigma}$ of $A$), then the
feasible regions of the linear programs in $LP_{A,c}$ have integral
vertices for each $b \in cone(A) \cap \mathbb Z^d$, and $yA
\leq c$ is TDI {\em for all} $c$. Hence if $A$ is unimodular, then
$IP_{A,c}$ is a Gomory family for all generic cost vectors
$c$. However, just as integrality of the optimal solutions of programs
in $LP_{A,c}$ is not necessary for $IP_{A,c}$ to be a Gomory family,
unimodularity of $A$ is not necessary for $IP_{A,c}$ to be a Gomory
family {\em for all} $c$. 

\begin{example}\label{allinishaveembeds}
Consider the seven by twelve integer matrix 
$$ A = \left[ \begin{array}{cccccccccccc}
1&0&0&0&0&0&1&1&1&1&1&0 \\
0&1&0&0&0&0&1&1&0&0&0&1 \\
0&0&1&0&0&0&1&0&1&0&0&1 \\
0&0&0&1&0&0&0&1&0&1&0&0 \\
0&0&0&0&1&0&0&0&1&0&1&0 \\
0&0&0&0&0&1&0&0&0&1&1&1 \\
0&0&0&0&0&0&1&1&1&1&1&1 \\
\end{array} \right]$$ 
of rank seven. The maximal minors of $A$ have absolute values zero,
one and two and hence $A$ is not unimodular. This matrix has $376$
distinct regular triangulations supporting $418$ distinct order ideals
$\mathcal O_c$ (computed using TiGERS). In each case, the standard pairs
of $\mathcal O_c$ are indexed by just the maximal simplices of the
regular triangulation $\Delta_c$ that supports it. Hence $IP_{A,c}$ is
a Gomory family for all generic $c$.  \qed\\
\end{example}

The above discussion shows that $IP_{A,c}$ being a Gomory family is
more general than $yA \leq c$ being TDI. Similarly, $IP_{A,c}$ being a
Gomory family for all generic $c$ is more general than $A$ being a
unimodular matrix.

\section{Gomory Families and Hilbert Bases}

As we just saw, unimodular matrices or more generally, unimodular
regular triangulations lead to Gomory families. A common property of
unimodular matrices and matrices $A$ such that $cone(A)$ has a
unimodular triangulation is that the columns of $A$ form a {\em
  Hilbert basis} for $cone(A)$, i.e., $\mathbb{N}A = cone(A) \cap
\mathbb{Z}^d$ (assuming $\mathbb Z A = \mathbb Z^d$).

\begin{definition} \label{normal}
A $d \times n$ integer matrix $A$ is {\em normal} if the semigroup
$\mathbb{N}A$ equals $cone(A) \cap \mathbb{Z}^d$.
\end{definition}

The reason for this (highly over used) terminology here is that if the
columns of $A$ form a Hilbert basis, then the zero set of the toric
ideal $I_A$ (called a {\em toric variety}) is a {\em normal} variety.
See \cite[Chapter 14]{Stu} for more details. We first note that if $A$
is not normal, then $IP_{A,c}$ need not be a Gomory family for any
cost vector $c$.

\begin{example} \label{nonnormal}
  The matrix $A = \left[\begin{array}{cccc} 1&1&1&1\\0&1&3&4
    \end{array}\right]$ is not normal since $(1,2)^t$ which lies in
  $cone(A) \cap \mathbb Z^2$ cannot be written as a non-negative
  integer combination of the columns of $A$. This matrix gives rise to
  10 distinct order ideals $\mathcal O_c$ supported on its four
  regular triangulations $\{\{1,4\}\},
  \{\{1,2\},\{2,4\}\},\{\{1,3\},\{3,4\}\}$ and
  $\{\{1,2\},\{2,3\},\{3,4\}\}$. Each $\mathcal O_c$ has at least one
  standard pair that is indexed by a lower dimensional face of
  $\Delta_c$.
\end{example} 

The matrix in Example~\ref{long-chain} is also not normal and has no
Gomory families. While we do not know whether normality of $A$ is
sufficient for the existence of a generic cost vector $c$ such that
$IP_{A,c}$ is a Gomory family, we will now show that under certain
additional conditions, normal matrices do give rise to Gomory
families.

\begin{definition} A $d \times n$ integer matrix $A$ is
  $\Delta$-normal if $cone(A)$ has a triangulation $\Delta$ such that
  for every maximal face $\sigma \in \Delta$, the columns of $A$ in
  $cone(A_{\sigma})$ form a Hilbert basis.
\end{definition}

\begin{remark}
  If $A$ is $\Delta$-normal for some triangulation $\Delta$, then it
  is normal. To see this note that every lattice point in $cone(A)$
  lies in $cone(A_{\sigma})$ for some maximal face $\sigma \in
  \Delta$. Since $A$ is $\Delta$-normal, this lattice point also lies
  in the semigroup generated by the columns of $A$ in $cone(A_\sigma)$
  and hence in $\mathbb N A$.
 
  Observe that $A$ is $\Delta$-normal with respect to all the
  unimodular triangulations of $cone(A)$. Hence triangulations
  $\Delta$ with respect to which $A$ is $\Delta$-normal generalize
  unimodular triangulations of $cone(A)$.
\end{remark}

Examples \ref{first} and \ref{second} show that the set of 
matrices where $cone(A)$ has a unimodular triangulation 
is a proper subset of the set of $\Delta$-normal matrices which in
turn is a proper subset of the set of normal matrices.

\begin{example} \label{first}
  Examples of normal matrices with no unimodular triangulations can be
  found in \cite{BoGo} and \cite{FZ}. If $cone(A)$ is simplicial for
  such a matrix, $A$ will be $\Delta$-normal with respect to its
  coarsest (regular) triangulation $\Delta$ consisting of the single
  maximal face with support $cone(A)$.  For
  instance, consider the following example taken from \cite{FZ}:
$$A = \left[ \begin{array}{cccccccc} 1 & 0 & 0 & 1 & 1 & 1 & 1 & 1
\\0 & 1 & 0 & 1 & 1 & 2 & 2 & 2 \\0 & 0 & 1 & 1 & 2 & 2 & 3
& 3 \\0 & 0 & 0 & 1 & 2 & 3 & 4 & 5 
\end{array} \right].$$
Here $cone(A)$ has $77$ regular triangulations and no
unimodular triangulations. Since $cone(A)$ is simplicial, $A$ is
$\Delta$-normal with respect to its coarsest regular triangulation
$\{\{1,2,3,8\}\}$.
\end{example}

\begin{example} \label{second}
There are normal matrices $A$ that are not $\Delta$-normal with respect to
any triangulation of $cone(A)$. To see such an example, consider
the following modification of the matrix in Example~\ref{first} 
that appears in \cite[Example 13.17]{Stu} :
$$A = \left[ \begin{array}{ccccccccc} 0 & 1 & 0 & 0 & 1 & 1 & 1 & 1 & 1
\\ 0 & 0 & 1 & 0 & 1 & 1 & 2 & 2 & 2 \\ 0 & 0 & 0 & 1 & 1 & 2 & 2 & 3
& 3 \\ 0 & 0 & 0 & 0 & 1 & 2 & 3 & 4 & 5 \\ 1 & 1 & 1 & 1 & 1 & 1 & 1
& 1 & 1
\end{array} \right].$$
This matrix is again normal and each of its nine columns generate an 
extreme ray of $cone(A)$. Hence the only way for this matrix to be
$\Delta$-normal for some $\Delta$ would be if $\Delta$ is a 
unimodular triangulation of $cone(A)$. However, there are no
unimodular triangulations in this example.
\end{example}

\begin{theorem} \label{specialinitial}
If $A$ is $\Delta$-normal for some regular triangulation $\Delta$ then
there exists a generic cost vector $c \in \mathbb Z^n$ such that
$\Delta = \Delta_c$ and $IP_{A,c}$ is a Gomory family.
\end{theorem}

\begin{proof}
Without loss of generality we can assume that the columns of $A$ in
$cone(A_{\sigma})$ form a minimal Hilbert basis for every 
maximal face $\sigma$ of $\Delta$. If there were a
redundant element, the smaller matrix obtained by removing this
column from $A$ would still be $\Delta$-normal.

For a maximal face $\sigma \in \Delta$, let $\sigma_{in} \subset
\{1,\ldots,n\}$ be the set of indices of all columns of $A$ lying in
$cone(A_{\sigma})$ that are different from the columns of
$A_{\sigma}$. Suppose $a_{i_1}, \ldots, a_{i_k}$ are the columns of
$A$ that generate the one dimensional faces of $\Delta$, and $c' \in
\mathbb R^n$ a cost vector such that $\Delta = \Delta_{c'}$. We modify
$c'$ to obtain a new cost vector $c \in \mathbb R^n$ such that $\Delta
= \Delta_c$ as follows. For $j = 1, \ldots, k$, let $c_{i_j} :=
c'_{i_j}$. If $j \in \sigma_{in}$ for some maximal face $\sigma \in
\Delta$, then $a_j =\sum_{i \in \sigma} \lambda_ia_i$, $0 \leq
\lambda_i < 1$ and we define $c_j := \sum_{i \in \sigma}
\lambda_ic_i$. Hence, for all $j \in \sigma_{in}$, $(a_j^t, c_j) \in
\mathbb R^{d+1}$ lies in $C_{\sigma} := cone((a_i^t,c_i) : i \in \sigma)
= cone((a_i^t,{c'}_i) : i \in \sigma)$ which was a facet of $C =
cone((a_i^t,{c'}_i) : i = 1, \ldots, n)$. If $y \in \mathbb R^d$ is a
vector as in Definition~\ref{regtriang} showing that $\sigma$ is a
maximal face of $\Delta_{c'}$ then $y \cdot a_i = c_i$ for all $i \in
\sigma \cup \sigma_{in}$ and $y \cdot a_j < c_j$ otherwise. Since
$cone(A_{\sigma}) = cone(A_{\sigma \cup \sigma_{in}})$, we conclude
that $cone(A_{\sigma})$ is a maximal face of $\Delta_c$.

If $b \in \mathbb{N}A$ lies in $cone(A_\sigma)$ for a maximal face
$\sigma \in \Delta_c$, then $IP_{A,c}(b)$ has at least one feasible
solution $u$ with support in $\sigma \cup \sigma_{in}$ since $A$ is
$\Delta$-normal. Further, $(b^t, c \cdot u) = ((Au)^t,c \cdot u)$ lies in
$C_{\sigma}$ and all feasible solutions of $IP_{A,c}(b)$ with support
in $\sigma \cup \sigma_{in}$ have the same cost value by
construction. Suppose $v \in \mathbb{N}^n$ is any feasible solution of
$IP_{A,c}(b)$ with support not in $\sigma \cup \sigma_{in}$. Then $c
\cdot u < c \cdot v$ since $(a_i^t,c_i) \in C_{\sigma}$ if and only if
$i \in \sigma \cup \sigma_{in}$ and $C_{\sigma}$ is a lower facet 
of $C$. Hence the optimal solutions of
$IP_{A,c}(b)$ are precisely those feasible solutions with support in
$\sigma \cup \sigma_{in}$. The vector $b$ can be expressed as $b = b'
+ \sum_{i \in \sigma} z_ia_i$ where $z_i \in \mathbb N $ are unique
and $b' \in \{\sum_{i \in \sigma} \lambda_i a_i \, : \, 0 \leq
\lambda_i < 1\} \cap \mathbb{Z}^d$ is also unique. The vector $b' =
\sum_{j \in \sigma_{in}} r_ja_j$ where $r_j \in \mathbb N$. Setting
$u_i = z_i$ for all $i \in \sigma$, $u_j = r_j$ for all $j \in
\sigma_{in}$ and $u_k = 0$ otherwise, we obtain all feasible solutions
$u$ of $IP_{A,c}(b)$ with support in $\sigma \cup \sigma_{in}$.

If there is more than one such feasible solution, then $c$ is not
generic. In this case, we can perturb $c$ to a generic cost vector
$c'' = c + \epsilon \omega$ by choosing $1 \gg \epsilon > 0$, $\omega_j
\ll 0$ whenever $j = i_1, \ldots, i_k$ and $\omega_j = 0$ otherwise.
Suppose $u_1, \ldots, u_t$ are the optimal solutions of the integer
programs $IP_{A,c''}(b')$ where $b' \in \{\sum_{i \in \sigma}
\lambda_i a_i \, : \, 0 \leq \lambda_i < 1\} \cap \mathbb{Z}^d$. (Note
that $t = |\{\sum_{i \in \sigma} \lambda_i a_i \, : \, 0 \leq
\lambda_i < 1\} \cap \mathbb{Z}^d|$ is the index of $\mathbb Z
A_{\sigma}$ in $\mathbb Z A$.)  The support of each such $u_i$ is
contained in $\sigma_{in}$.  For any $b \in cone(A_{\sigma}) \cap
\mathbb Z^d$, the optimal solution of $IP_{A,c''}(b)$ is hence $u =
u_i + z$ for some $i \in \{1, \ldots, t\}$ and $z \in \mathbb
N^n$ with support in $\sigma$. This shows that $\mathbb N A$
is covered by the affine semigroups $\phi_A(S(u_i,\sigma))$ where
$\sigma$ is a maximal face of $\Delta$ and $u_i$ as above for each
$\sigma$. By construction, the corresponding admissible pairs
$(u_i,\sigma)$ are all standard for $\mathcal O_{c''}$. Since all data 
is integral, $c'' \in \mathbb Q^n$ and hence can be scaled to lie 
in $\mathbb Z^n$. Renaming $c''$ as $c$, we conclude that
$IP_{A,c}$ is a Gomory family. 
\end{proof}

\begin{corollary} Let $A$ be a normal matrix such that $cone(A)$ is
  simplicial, and let $\Delta$ be the coarsest triangulation 
  whose single maximal face has support $cone(A)$. Then there exists a
  cost vector $c \in \mathbb Z^n$ such that $\Delta = \Delta_c$ and 
  $IP_{A,c}$ is a Gomory family.  
\end{corollary}

\begin{example} \label{non-delta-Gomory}
Consider the normal matrix in Example \ref{example-with-gfamily}.
Here $cone(A)$ is generated by the first, second and sixth columns of
$A$ and hence $A$ is $\Delta$-normal with respect to the regular
triangulation $\{\{1, \, 2, \, 6\}\}$. There are 13 distinct sets
$\mathcal O_c$ supported on $\Delta$. Among the 13 corresponding
families of integer programs, only one is a Gomory family. A
representative cost vector for this $IP_{A,c}$ is
$c=(0,0,4,4,1,0)$. The standard pair decomposition of $\mathcal O_c$
is the one constructed in Theorem \ref{specialinitial}. The affine
semigroups $S(\cdot, \sigma)$ from this decomposition are:
$$ S(0, \sigma), \,\, S(e_3, \sigma), \,\, S(e_4, \sigma), \,\,
\text{and} \,\, S(e_5, \sigma).$$ 
Note that $A$ is not $\Delta$-normal
with respect to the regular triangulation supporting the Gomory family
$IP_{A,c}$ in Example \ref{example-with-gfamily}. The columns of $A$
in $cone(A_{\sigma_1})$ are the columns of $A_{\sigma_1}$ and $A_3$.
The vector $(1, \, 2, \, 2)$ is in the minimal Hilbert basis of 
$cone(A_{\sigma_1})$ but is not a column of $A$. This example shows
that a regular triangulation $\Delta$ of $cone(A)$ can support a
Gomory family even if $A$ is not $\Delta$-normal. The Gomory families 
in Theorem~\ref{specialinitial} have a very special standard pair 
decomposition. \qed\\
\end{example}

\begin{problem} \label{Gomory-question}
If $A \in \mathbb Z^{d \times n}$ is a normal matrix, does there
exist a generic cost vector $c \in \mathbb Z^n$ such that $IP_{A,c}$
is a Gomory family? 
\end{problem}

While we do not know the answer to this question, we will now show
that stronger results are possible for small values of $d$. 

\begin{theorem} \label{unimodular-for-three}
If $A \in \mathbb Z^{d \times n}$ is a normal matrix and $d \leq 3$,
then there exists a generic cost vector $c \in \mathbb Z^n$ such that
$IP_{A,c}$ is a Gomory family.
\end{theorem}

\begin{proof}
It is known that if $d \leq 3$ then $cone(A)$ has a regular unimodular 
triangulation $\Delta_c$ \cite{Seb}. The result then follows from 
Corollary \ref{TDI-Gomory}.
\end{proof}

Before we proceed, we rephrase Problem~\ref{Gomory-question} in terms
of covering properties of $cone(A)$ and $\mathbb N A$ along the lines
of \cite{BoGo}, \cite{BG}, \cite{BGHMW}, \cite{FZ} and \cite{Seb}. To
obtain the same set up as in these papers we assume in this section
that $A$ is normal and the columns of $A$ form the unique minimal
Hilbert basis of $cone(A)$. Using the terminology in \cite{BG}, the
{\em free Hilbert cover} problem asks whether there exists a covering
of $\mathbb N A$ by semigroups $\mathbb{N}A_\tau$ where the columns of
$A_\tau$ are linearly independent. The {\em unimodular Hilbert cover}
problem asks whether $cone(A)$ can be covered by full 
dimensional unimodular subcones
$cone(A_{\tau})$ (i.e., $\mathbb Z A_{\tau} = \mathbb Z^d$), while the
stronger {\em unimodular Hilbert partition} problem asks whether
$cone(A)$ has a unimodular triangulation. (Note that if $cone(A)$ has
a unimodular Hilbert cover or partition using subcones
$cone(A_{\tau})$, then $\mathbb N A$ is covered by the semigroups
$\mathbb N A_{\tau}$.)  All these problems have positive answers if $d
\leq 3$ since $cone(A)$ admits a unimodular Hilbert partition in this
case \cite{BoGo}, \cite{Seb}. Normal matrices (with $d = 4$) such that
$cone(A)$ has no unimodular Hilbert partition can be found in
\cite{BoGo} and \cite{FZ}. Examples (with $d = 6$) that admit no free
Hilbert cover and hence no unimodular Hilbert cover can be found in
\cite{BG} and \cite{BGHMW}.

When $yA \leq c$ is TDI, the standard pair decomposition of $\mathbb N
A$ induced by $c$ gives a unimodular Hilbert partition of $cone(A)$ by
Theorem~\ref{TDI-ness}. An important difference between
Problem~\ref{Gomory-question} and the Hilbert cover problems is that
{\em affine} semigroups cannot be used in Hilbert covers. Moreover,
affine semigroups that are allowed in standard pair decompositions
come from integer programming. If there are no restrictions on the
affine semigroups that can be used in a cover, $\mathbb N A$ can
always be covered by full dimensional affine semigroups: for any
triangulation $\Delta$ of $cone(A)$ with maximal subcones
$cone(A_{\sigma})$, the affine semigroups $b + \mathbb{N}A_\sigma$
cover $\mathbb N A$ as $b$ varies in $\{ \sum_{i \in \sigma} \lambda_i
a_i: \,\, 0 \leq \lambda_i < 1\} \cap \mathbb{Z}^d$ and $\sigma$
varies among the maximal faces of the triangulation. A {\em partition}
of $\mathbb N A$ derived from this idea can be found in \cite[Theorem
5.2]{Sta82}. We recall the notion of {\em supernormality} 
introduced in \cite{HMS}.

\begin{definition} \label{supernormal}
A matrix $A \in \mathbb Z^{d \times n}$ is {\em supernormal}
if for every submatrix $A'$ of $A$, the columns of $A$ that 
lie in $cone(A')$ form a Hilbert basis for $cone(A')$.  
\end{definition}

\begin{proposition} \label{equivalence}
For $A \in \mathbb Z^{d \times n}$, the following are equivalent:
\begin{enumerate}
\item[(i)] $A$ is supernormal,
\item[(ii)] $A$ is $\Delta$-normal for every regular triangulation
$\Delta$ of $cone(A)$,
\item[(iii)] Every triangulation of $cone(A)$ in which all columns of 
$A$ generate one dimensional faces is unimodular.
\end{enumerate}
\end{proposition}

\begin{proof} The equivalence of (i) and (iii) was established in 
\cite[Proposition 3.1]{HMS}. Definition \ref{supernormal} shows that
(i) $\Rightarrow$ (ii). Hence we just need to show that (ii)
$\Rightarrow$ (i). Suppose that $A$ is $\Delta$-normal for every
regular triangulation of $cone(A)$. In order to show that $A$ is supernormal
we only need to check submatrices $A'$ where the dimension of
$cone(A')$ is $d$.  Choose a cost vector $c$ with $c_i \gg 0$ if the
$i$th column of $A$ does not generate an extreme ray of $cone(A')$,
and $c_i = 0$ otherwise. This gives a polyhedral subdivision of
$cone(A)$ in which $cone(A')$ is a maximal face. There are
standard procedures that will refine this subdivision to a regular
triangulation $\Delta$ of $cone(A)$. Let $T$ be the set of maximal faces
$\sigma$ of $\Delta$ such that $cone(A_{\sigma})$ lies in
$cone(A')$. Since $A$ is $\Delta$-normal, the columns of $A$ that lie
in $cone(A_{\sigma})$ form a Hilbert basis for $cone(A_{\sigma})$ for
each $\sigma \in T$. However, since their union is the set of columns
of $A$ that lie in $cone(A')$, this union forms a Hilbert basis for
$cone(A')$.
\end{proof}

It is easy to catalog all $\Delta$-normal and supernormal matrices, of
the type considered in this paper, for small values of $d$. We say
that the matrix $A$ is {\em graded} if its columns span an affine
hyperplane in $\mathbb R^d$. If $d=1$, $cone(A)$ has $n$
triangulations $\{\{i\}\}$ each of which has the unique maximal
subcone $cone(a_i)$ whose support is $cone(A)$. If 
we assume that $a_1 \leq a_2 \leq \cdots \leq a_n$, then $A$ is normal
if and only if either $a_1 = 1$, or $a_n = -1$. Also, $A$ is normal if
and only if it is supernormal. If $d=2$ and the columns of $A$ are
ordered counterclockwise around the origin, then $A$ is normal if and
only if $det(a_i, a_{i+1}) = 1$ for all $i = 1, \ldots, n-1$. Such an
$A$ is supernormal since it is $\Delta$-normal for every triangulation
$\Delta$ --- the Hilbert basis of a maximal subcone of $\Delta$ is
precisely the set of columns of $A$ in that subcone. If $d = 3$ then
as mentioned before, $cone(A)$ has a unimodular triangulation with
respect to which $A$ is $\Delta$-normal. However, not every such $A$
needs to be supernormal: the matrix in Example
\ref{example-with-gfamily} is not $\Delta$-normal for the $\Delta$ 
supporting the Gomory family in that example.
If $d=3$ and $A$ is graded, then
without loss of generality we can assume that the columns of $A$ span
the hyperplane $x_1 = 1$. If $A$ is normal as well, then its columns
are precisely all the lattice points in the convex hull of
$A$. Conversely, every graded normal $A$ with $d=3$ arises this way ---
its columns are all the lattice points in a polygon in $\mathbb R^2$
with integer vertices. In particular, every triangulation of $cone(A)$
that uses all the columns of $A$ is unimodular. Hence, by Proposition
\ref{equivalence}, $A$ is supernormal, and therefore $\Delta$-normal
for any triangulation of $A$. 

\begin{theorem}\label{smalld}
Let $A \in \mathbb Z^{d \times n}$ be a normal matrix of rank $d$.

\begin{enumerate}

\item[(i)] If $d = 1,2$ or $A$ is graded and $d = 3$, every
  regular triangulation of $cone(A)$ supports at least one
  Gomory family.

\item[(ii)] If $d = 2$ and $A$ is graded, every regular
  triangulation of $cone(A)$ supports exactly one Gomory family.

\item [(iii)] If $d = 3$ and $A$ is not graded, or if 
  $d = 4$ and $A$ is graded, then not all regular
  triangulations of $cone(A)$ may support a Gomory family. In
  particular, $A$ may not be $\Delta$-normal with respect to every
  regular triangulation. 
\end{enumerate}
\end{theorem}

\begin{proof}
(i) If $d=1,2$ or $A$ is graded and $d=3$, $A$ is supernormal and
hence by Proposition~\ref{equivalence} and
Theorem~\ref{specialinitial}, every regular triangulation of $cone(A)$
supports at least one Gomory family.

\smallskip
\noindent
(ii) If $d=2$ and $A$ is graded, then we may assume that
$$A = \left[ \begin{array}{ccccc}
    1 & 1 & 1 & \ldots & 1 \\
    0 & 1 & 2 & \ldots & n-1 \end{array} \right].$$
In this case, $A$
is supernormal and hence every regular triangulation $\Delta$ of
$cone(A)$ supports a Gomory family by Theorem~\ref{specialinitial}.
Suppose the maximal cones of $\Delta$, in counter-clockwise order, are
$C_1, \ldots, C_r$. Assume the columns of $A$ are labeled
such that $C_i = cone(a_{i-1},a_i)$ for $i = 1, \ldots, r$, and 
the columns of $A$ in the interior of $C_i$ are labeled in
counter-clockwise order as $b_{i1}, \ldots, b_{ik_i}$. Hence the $n$
columns of $A$ from left to right are:
$$a_0, b_{11}, \cdots, b_{1k_1}, a_{1}, b_{21}, \cdots, a_{r-1},
b_{r1}, \cdots, b_{rk_r}, a_{r}.$$
Indexing the columns of $A$ by
their labels, the maximal faces of $\Delta$ are $\sigma_i = \{
i-1,i\}$ for $i = 1, \ldots, r$. Let $e_i$ be the unit vector of
$\mathbb R^n$ indexed by the true column index of $a_i$ in $A$ and
$e_{ij}$ be the unit vector of $\mathbb R^n$ indexed by the true
column index of $b_{ij}$ in $A$.  Since the columns of $A$ form a
minimal Hilbert basis of $cone(A)$, $e_i$ is the unique solution to
$IP_{A,c}(a_i)$ for all $c$ and $e_{ij}$ is the unique solution to
$IP_{A,c}(b_{ij})$ for all $c$.  Hence the standard pairs of
Theorem~\ref{specialinitial} are $(0,\sigma_i)$ and
$(e_{ij},\sigma_i)$ for $i=1, \ldots, r$ and $j = 1, \ldots, k_i$.

Suppose $\Delta$ supports a second Gomory family $IP_{A,\omega}$.
Then every standard pair of ${\mathcal O}_w$ is also of the form
$(\ast, \sigma_i)$ for $\sigma_i \in \Delta$, and $r$ of them are
$(0, \sigma_i)$ for $i = 1, \ldots, r$. The remaining standard pairs
are of the form $(e_{ij}, \sigma_k)$. To see this, consider the
semigroups in $\mathbb N A$ arising from the standard pairs of
$\mathcal O_w$. The total number of standard pairs of $\mathcal
O_c$ and $\mathcal O_w$ are the same. Since the columns of $A$ all lie 
on $x_1 =1$, no two $b_{ij}$s can be covered by a semigroup
coming from the same standard pair and none of them are covered by a
semigroup $(0, \sigma_i)$. We show that if $(e_{ij},
\sigma_k)$ is a standard pair of $\mathcal O_w$ then $k = i$ and thus
$\mathcal O_w = \mathcal O_c$.

If $r=1$, the standard pairs of $\mathcal O_w$ are $(0,\sigma_1),
(e_{11}, \sigma_1), \ldots, (e_{1k_1}, \sigma_1)$ as in
Theorem~\ref{specialinitial}. If $r > 1$, consider the last cone $C_r
= cone(a_{r-1}, a_r)$. If $a_{r-1}$ is the second to last column of $A$,
then $C_r$ is unimodular and the semigroup from $(0,\sigma_r)$ covers
$C_r \cap \mathbb Z^2$. The subcomplex comprised of $C_1, \ldots,
C_{r-1}$ is a regular triangulation $\Delta'$ of $cone(A')$ where $A'$
is obtained by dropping the last column of $A$. Since $A'$ is a normal
graded matrix with $d=2$ and $\Delta'$ has less than $r$ maximal
cones, the standard pairs supported on $\Delta'$ are as in
Theorem~\ref{specialinitial} by induction. If $a_{r-1}$ is not the
second to last column of $A$ then $b_{rk_r}$, the second to last column of
$A$ is in the Hilbert basis of $C_r$ but is not a generator of $C_r$.
So ${\mathcal O}_w$ has a standard pair of the form $(e_{rk_r},
\sigma_i)$.  If $\sigma_i \neq \sigma_r$, then the lattice point
$b_{rk_r} + a_r$ cannot be covered by the semigroup from this or any
other standard pair of $\mathcal O_w$. Hence $\sigma_i = \sigma_r$. By
a similar argument, the remaining standard pairs indexed by $\sigma_r$
are $(e_{r (k_r-1)}, \sigma_r), \ldots,(e_{r1},\sigma_r)$ along with
$(0,\sigma_r)$. These are precisely the standard pairs of $\mathcal
O_c$ indexed by $\sigma_r$. Again we are reduced to considering the
subcomplex comprised of $C_1, \ldots, C_{r-1}$ and by induction, the
remaining standard pairs of $\mathcal O_w$ are as in
Theorem~\ref{specialinitial}.

\smallskip
\noindent
(iii) The 
$3 \times 6$ normal matrix $A$  of Example \ref{example-with-gfamily} 
has 10 distinct Gomory families
supported on 10 out of the 14 regular triangulations of $cone(A)$.
Furthermore, the normal matrix $$A = \left[ \begin{array}{ccccccc}
1& 1& 1& 1& 1& 1& 1 \\
1& 0& 1& 1& 1& 1& 0\\
0& 1& 2& 2& 1& 1& 0\\
0& 0& 4& 3& 2& 1& 0 \end{array} \right]$$ has 11 distinct Gomory 
families supported on 11 out of its 19 regular triangulations.
\end{proof}

\bibliography{refs} \bibliographystyle{plain}

\end{document}